# 数控加工中路径规划与速度插补综述

马鸿宇[1]，申立勇[1]，姜　鑫[2]，邹　强[3]，袁春明[4]

(1. 中国科学院大学数学科学学院，北京 100049；2. 北京航空航天大学人工智能研究院，数学、信息与行为教育部重点实验室，北京 100083；3. 浙江大学计算机科学与技术学院，CAD&CG 国家重点实验室，浙江 杭州 310027；4. 中国科学院数学与系统科学研究院，北京 100190)

摘　　要：数控技术在现代制造工业中被广泛使用，相关研究一直为学界和业界共同关注。数控技术的传统流程主要包含刀具路径规划和进给速度插补。为实现高速高精加工，人们通常将路径规划与速度插补中的若干问题转换成数理优化模型，针对工程应用问题的复杂性，采用分步迭代优化的思路进行求解，但所得的结果往往只是局部最优解。其次，路径规划与速度插补都是为了加工一个工件曲面，分两步进行处理虽然简化了计算，但也导致不能进行整体优化。因此，为了更好地开展路径规划与速度插补一体化设计与全局最优求解的研究，系统性地了解并学习已有的代表性工作是十分有必要的。所以将逐次介绍数控加工中刀具路径规划与速度插补的相关方法与技术进展，包括基于端铣的加工路径规划；刀轴方向优化；G 代码加工以及拐角过渡；参数曲线路径的进给速度规划等国内外相关研究以及最新提出的一些新型加工优化方法。

## 0 背景介绍

数控机床(computer numerical control machine tools)是数字化控制机床的简称，是一种装有程序控制系统的零件加工生产设备，可以高效生产各种复杂曲面零部件，如大型舰船螺旋桨、航空发动机叶轮、隐形战机一体成型机翼，汽车发动机零部件，5G 高精度零部件等。数控机床作为工业母机是支撑国民经济发展的战略性产业。数控机床加工技术具有加工精度高、自动化程度高、可靠性高、对加工对象的适应性强等优点，已经成为当前制造业中主要力量，也是未来先进制造加工领域中必不可少的核心组成。数控系统作为数控机床的控制核心，是数控机床的“大脑”，其功能、控制精度和可靠性直接影响机床的整体性能、性价比和市场竞争力。数控技术则是数控系统中技术概括，传统意义上涵盖数理模型、加工工艺、系统控制、软硬件协同等，在当前智能框架下还可以延伸到前端的计算机辅助设计(computer aided design，CAD)几何设计分析，后端的数控加工(computer numerical control，CNC)智能控制等。

随着国际竞争日益激烈，各国对高、精、尖装备的需求日益扩大，向现代制造加工技术提出了新的挑战，所以数控加工技术已经是各主要国家的重要产业战略布局点。美国先后发布了《先进制造伙伴关系计划》、《先进制造业国家战略计划》、《美国先进制造业领导战略》等多项法案计划来大力发展先进制造加工技术；德国提出了工业 4.0 的发展计划，力图将其打造为新一轮工业革命的技术平台；日本公布了产业结构蓝图，旨在强化国内制造业；英国、韩国、印度、中国台湾等地区亦提出积极的战略和政策，推动新兴技术在数控机床等装备产业领域加快融合。在 2015 年，由国务院印发的《中国制造 2025》中提到，要将高档数控机床和基础制造装备列为“加快突破的战略必争领域”，锚定我国装备制造业全球竞争地位，支撑国防和产业安全的战略需求，特别是支撑我国从制造大国到制造强国的过渡中，核心战略技术研发是重中之重。

数控加工的主要流程包括 CAD、计算机辅助制造(computer aided manufacturing，CAM)、后处理过程以及机床加工，如图 1 所示。

随着数控机床的快速发展，人们在实际数控加工中对 CAM 的精度和效率提出了更高的要求。CAM 一般包含加工路径规划和刀轴方向规划。由于刀具必须沿着规划好的加工轨迹行进，因此，加工路径和刀轴方向的优劣在一定程度上决定了后继的加工精度和加工效率。其中加工路径规划又可以分为基于侧铣和基于端铣的刀具路径。侧铣在加

工一些特定曲面(如直纹面等)具有很高地加工效率，而端铣则在加工自由曲面时具有更好的效果。本文将主要介绍针对自由曲面的相关数控加工技术，因此关于侧铣的方法将不做展开。

而机床加工过程的核心内容则是进给速度规划。根据刀具路径的类型可分为：G 代码加工及拐角过渡；针对参数曲线的进给速度规划。在得到刀具路径之后，总是希望刀具在满足运动学约束和几何误差约束的前提下能够更快地沿着刀具路径进行切削，从而缩短加工时间。所以进给速度规划直接决定了加工工件的质量与效率，在数控加工中起着至关重要的作用。

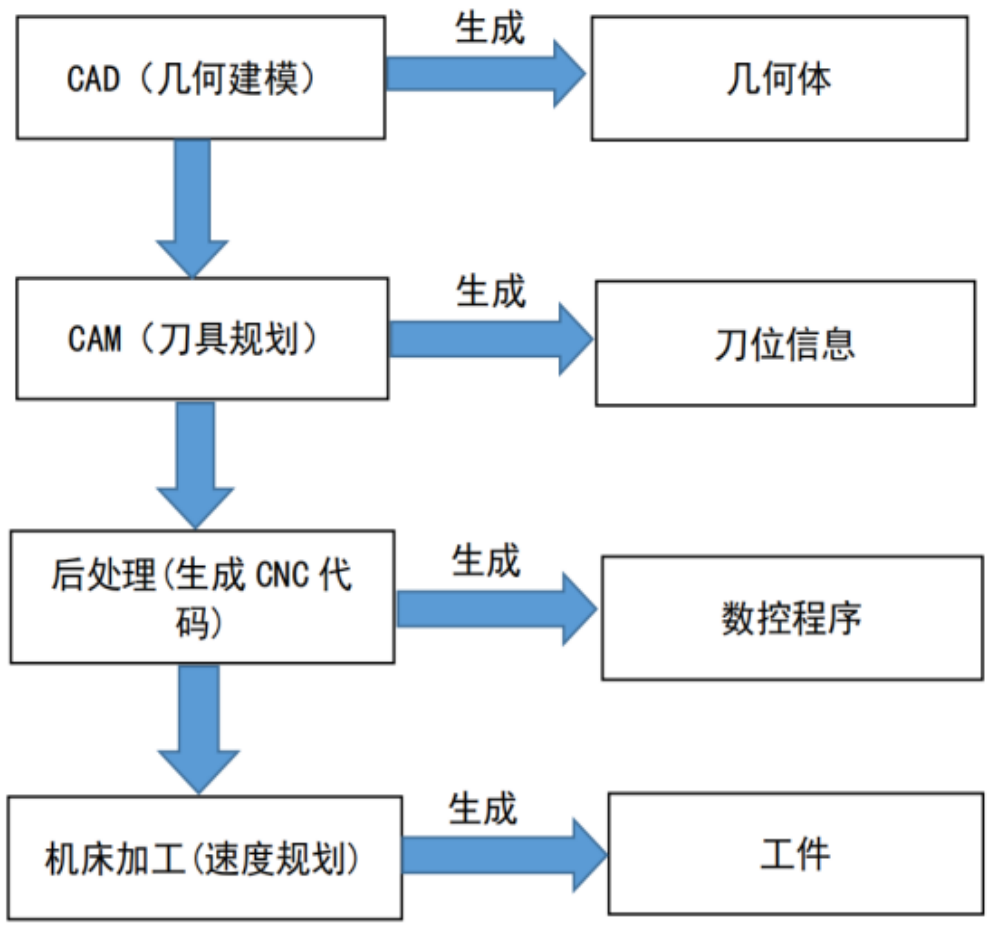

图 1　数控加工流程图

# 1　加工路径规划

加工路径指工件与机床刀具的相对运动轨迹，其决定着刀具如何切削毛坯，以形成设计几何所规定的曲面形状。加工路径在很大程度上决定着进给速度、切屑面积、切削力甚至是颤振(chatter)的发生与否，因此，路径质量直接关系着实际加工的精度和效率。从 CAD/CAM 一体化角度来看，加工路径是 CAD 模型(即设计几何)和 CNC 加工代码的中间状态(即制造几何)，是数字空间和物理空间之间的桥梁(图 2)。故而能否自动规划加工路径直接关系着产品生命周期的整体自动化程度和效率。

广义上，加工路径包括刀位点位置与刀轴方向 2 个方面；狭义上，加工路径仅指刀位点位置(一方面由于刀轴方向规划一般是在刀位点规划的基础之上进行，另一方面是因为 3 轴数控机床的广泛应用，这类机床的刀轴方向是固定不变的，不需要额外的刀轴方向规划)。因此，本文采用后者定义。加工路径的质量有 3 个主要评价标准[1]：

(1)　相邻路径间残留材料的高度(简称残高)不超过给定容差，以保证加工精度；

(2)　加工路径总长度尽量小，以减少加工时间；

(3)　加工路径的光滑程度尽量高，以防止频繁加减速，从而提高平均进给速度，降低加工时间。

其中，标准(1)是加工精度方面的要求，标准(2)和(3)是加工效率方面的要求。对于简单几何，3 个要求一般可以同时满足， 如对于平面，行型(direction parallel)的等残高路径同时也是最短且光顺的(除各条路径首尾相连处)；然而，对于复杂曲面，其之间往往是不相容的，不能同时满足。加工路径规划的难点就在于如何平衡三者，尤其是全局最优的平衡。

针对上述 3 个要求，人们在过去的 30 余年里提出了一系列路径规划方法，其中具有基础意义的加工路径生成方法主要有等参数法(iso-parametric)、截平面法(iso-planar)、等残高法(iso-scallop)等[2]，新近发展的加工路径优化方法有向量场方法、等水平法等。

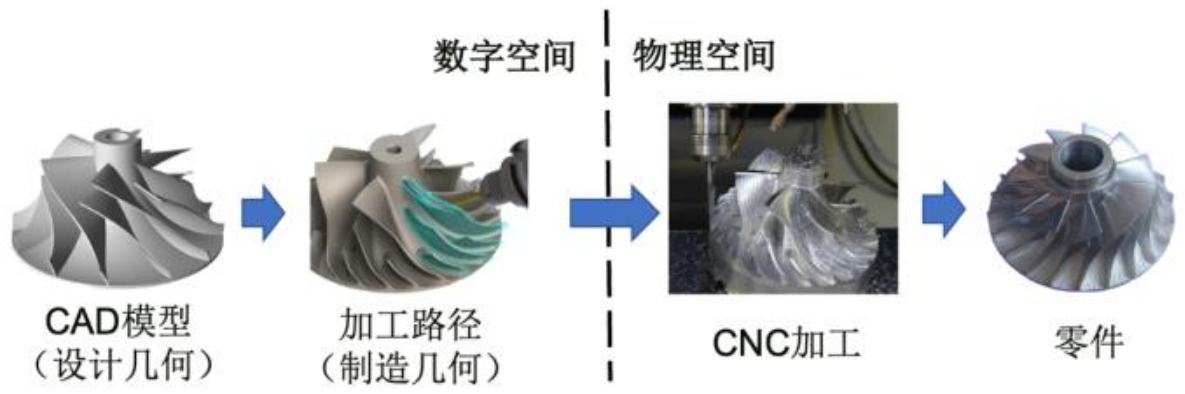

图 2　加工路径是数字-物理二元空间的桥梁

## 1.1　传统加工路径规划方法

在三大传统方法中，等参数法最早被学术界提出。如图 3 所示，对于三维空间中的自由曲面，该方法通过固定一个参数($u$)，变动另一参数($v$)，在参数域生成一条直线段，并选取其所对应的三维空间曲线为加工路径[3]。相邻加工路径对应的参数增量由工程师指定的容差确定，即残高不超过给定容差值。

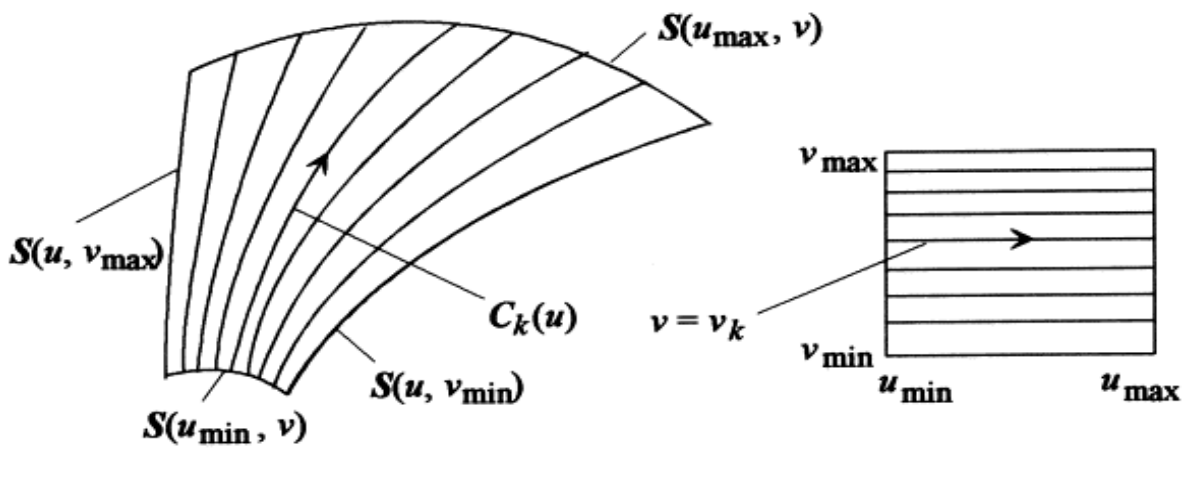

图 3　等参数路径规划方法示意图[3]

等参数法具有路径光滑、计算简单、与 NURBS

曲面表示方法契合度高等优点，但也存在路径质量不高的缺点，主要表现为路径总长度过大(其原因在于路径疏密不一致、冗余加工过多)，且该方法不适用于离散的网格和点云曲面。

针对前一缺点，常用的改进方法有：ELBER 和 COHEN[4]对过密区域进行路径裁剪；HE 等[5]对过密区域进行路径间距微调。针对后一缺点，可用网格模型和点云模型进行参数化，重建参数域解决；这方面典型的工作见文献[6-7]，其分别面向网格和点云。

通过上述改进方法，等参加工路径的总长度虽有所降低，但问题仍然突出。为进一步降低路径总长度，DING 等[8]提出截平面法[1]。其基本思路是将待加工自由曲面与一组平行面相交，然后以相交线为加工路径，如图 4 所示。除使用一组平面外，人们还尝试过使用一组曲面与待加工自由曲面进行相交以生成加工路径，如 KIM 和 CHOI [9]提出的引导面方法中的投影操作在几何上等价于求交。实际加工表明，截平面路径较等参路径具有更均匀的路径疏密度，更小的路径总长度。与此同时，截平面法具有可直接应用于网格和点云模型的优势，如文献[10-11]所提的方法。然而，截平面路径仍未彻底解决路径疏密不一致、冗余加工过多的问题。

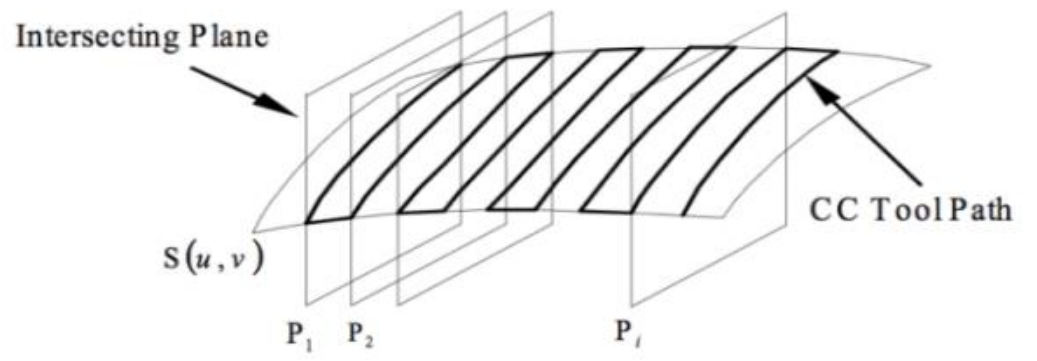

图 4　截平面路径规划方法示意图[8]

对上述问题，常用的改进思路有 2 个：曲面分割法和等残高法。其中，曲面分割法将待加工曲面分为多个区域，然后在不同区域生成具有不同疏密度的截平面加工路径，以此让各自区域内的路径疏密度尽量一致，其典型工作有 HU 等[12]提出的方法。显然，该方法只能在一定程度上减轻截平面法缺点带来的影响，不能完全解决问题，在各自区域内仍存在路径疏密不一致、冗余加工过多的问题。然而，曲面分割的思想在加工轨迹规划领域具有重要意义，最近这一思想又被重新挖掘出来，形成了近期基于方向场的加工轨迹规划方法[13]。

等残高法是一个能够完全解决路径疏密不一致和冗余加工问题的方法。其基本思路如图 5 所示，主要包括 2 部分：

(1)　选定一条初始路径(如某段边界)并以此为基础在待加工曲面上生成一簇偏移曲线；

(2)　在偏移时，相邻曲线间的残高保持恒定(即严格等于给定容差值，而不是像等参数法或截平面法中的不超过给定容差值)。

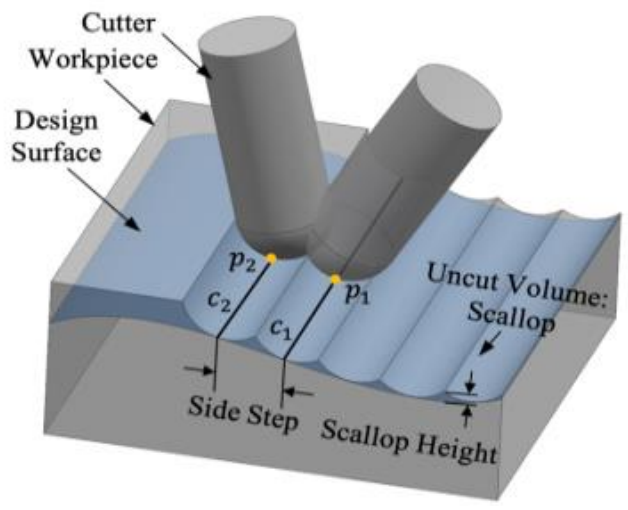

图 5　等残高路径规划方法示意图[19]

这一方法一经提出便受到广泛关注，多个改进方案随之被提出，FENG 和 LI[14]在提高计算精度方面做出工作，TOURNIER 和 DUC[15]在提高计算效率方面做了尝试，WEN 等[16]把工作应用范围方面做了拓广(如平底刀具，网格或点云曲面)。与此同时，针对不同的加工需求，人们还将之拓展到如图 6 所示的多种拓扑形式[17]，如螺旋形加工轨迹适用于高速加工[18]。

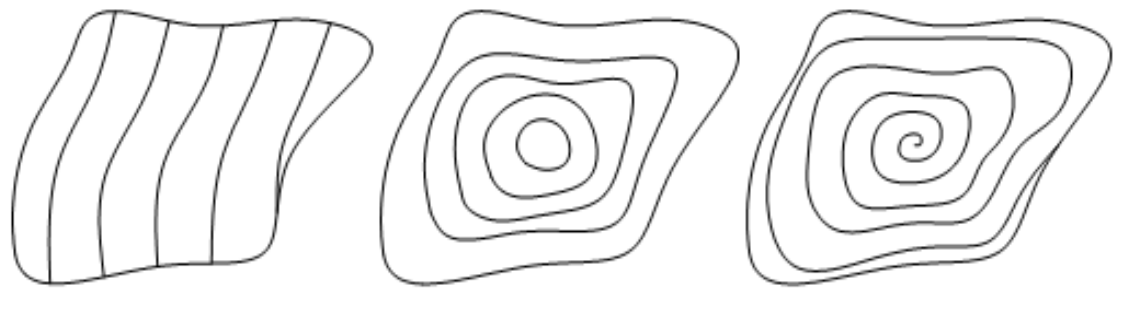

图 6　加工路径拓扑[7]((a)行型；(b)回型；(c)螺旋型)

等残高法因彻底解决了冗余加工问题，具有相对其他方法较少的加工路径总长度。然而，也存在两方面问题：

(1)　加工路径的总体形状依赖于初始路径的选取，其总长度也会有所不同，因此，还有改进空间；

(2)　等残高加工路径的生成依赖于偏移操作，而曲线偏移往往会在路径上产生尖角，如图 7 所示(蓝色为初始路径，黑色为偏移曲线)，因此等残高路径在光滑度方面还有优化空间。

[1] 其实在截平面法被学术界正式提出之前，工业界已独立应用该方法，但具体的算法步骤要么未披露，要么不够系统。

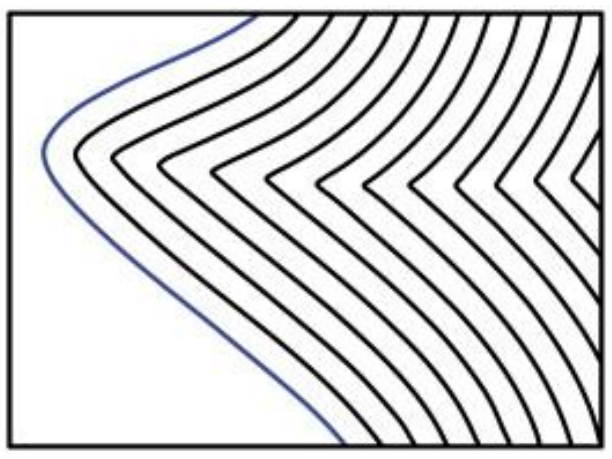

图 7　曲线偏移导致尖角问题

针对这 2 个问题，人们提出了一系列优化方法，如最近的基于方向场的加工路径生成方法[13]便是一种通过方向场的流线来对路径进行优化选择的方法。

## 1.2　新型加工路径优化方法

针对初始路径选择的问题,其使用一种类似贪心算法的策略——所选取的初始路径具有最优的加工条件，其他路径仍然使用传统偏移的方法得到。最优初始路径的选取标准也多种多样，经典的有最大加工带宽(图 8)、最大/最小曲面曲率等。这类贪心策略有一个明显的缺点，远离初始路径的加工路径的最优性难以保证。

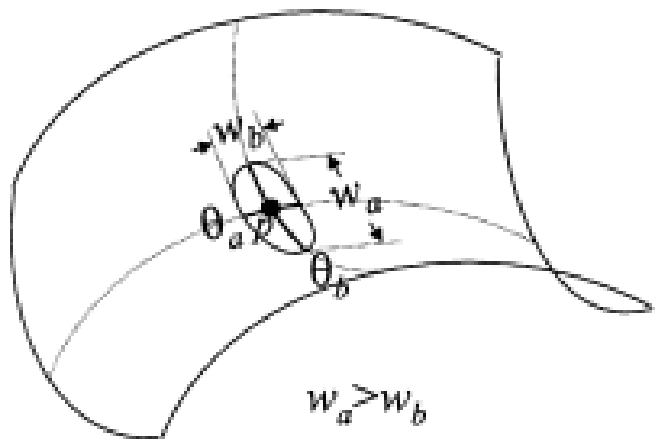

(a)

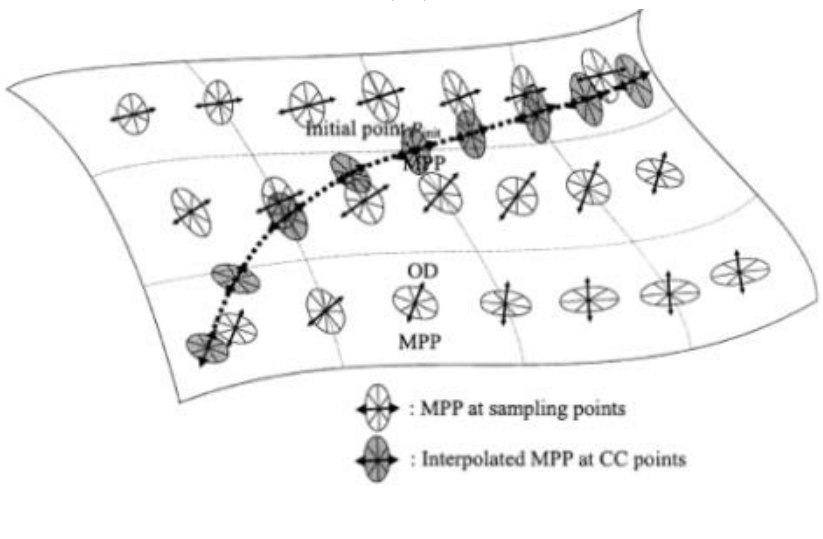

(b)

图 8　最大加工带宽示意图[2]((a)最大带宽方向示意图；(b)根据向量场生成最优初始路径)

上述问题的一大改进策略是对曲面进行分割，然后在各个被分割区域分别选取最优初始路径。由于被分割区域一般不大，初始路径与区域内的其他路径相距不远，其最优性在一定程度上能够被保证。传统上曲面分割方法是直接对曲率、加工带宽、运动学性能等加工性能指标进行自底向上的聚类或自顶向下的分类，最近的曲面分割的方法是通过将加工性能指标转化成方向场(很多文献中使用矢量场来描述这类方法，这是不严格确的，此类方法只利用了方向信息，不利用模长信息)，并利用场中的奇异点来对曲面进行分割，如图 9 和图 10 所示。场方法的依据在于奇异点分割出来的区域内部的流线一般具有很好的相似性[19]，这可以保持偏移路径与初始路径的一致性。

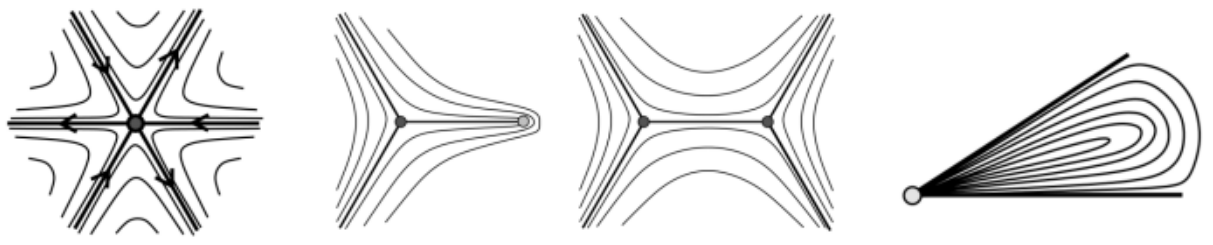

**图 9**　方向场奇异点示意图[19]

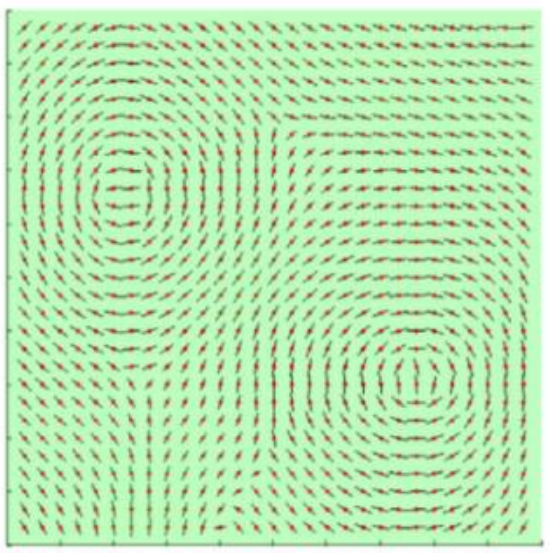

(a)

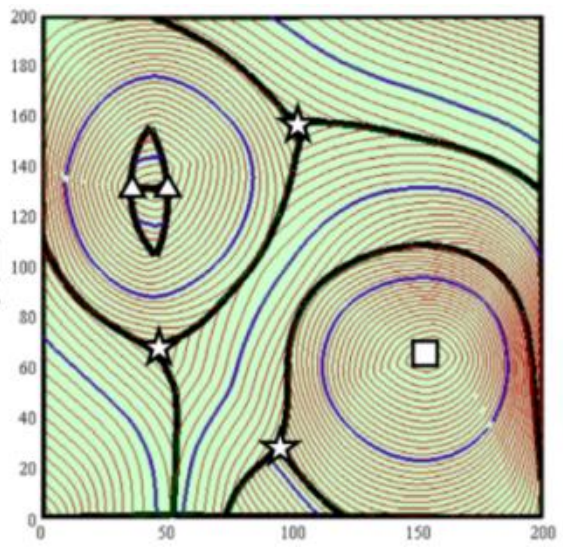

(b)

**图 10**　基于方向场的曲面分割方法[13]((a)代表最优加工方向的方向场；(b)根据奇异分界对曲面进行分割)

利用方向场(或者说张量场，本质一样)的奇异点来进行曲面分割这一思想首先由文献[13]提出。随后，学者提出一系列改进方法[20-24]，主要沿着以下 3 方面展开研究：

(1) 向量场生成方法：将生成方向场的加工性能指标从最大加工带宽推广至最大材料去除率[20]、最大进给速度[21]、最小能量消耗[22]；

(2) 分割区域内路径生成方法：将文献[13]中的初始路径人工选择方法扩展至自动选择方法[23]；

(3) 分割区域内路径拓扑确定方法：文献[24]将文献[13]中的人工指定路径拓扑方法扩展至根据流线形状自动确定路径拓扑的方法。

虽然基于方向场的加工路径规划方法取得了

很好的效果，方向场也具有表示多种加工指标的优点，但仍有一个重要问题有待回答：按奇异点来分割曲面是否是全局最优分割，分割区域内路径生成方法是否为全局最优？目前的方法在这两方面均采用启发式算法，其与全局最优路径的关系尚未知。最近 ZOU[25]提出了一种将方向场推广至矢量场的方法，首次实现了路径总长度的全局优化，但尚未包含其他加工指标的全局优化。总之，场方法未来仍有很大地改进空间。

与路径长度优化相比，有关路径光滑度的工作较少。初期的方法是在已有路径基础上对尖角部分进行局部修正[26]。这类方法具有实现简单的优点，但也存在 2 个问题：①在修正处破坏了等残高性质；②局部修正只能给出次优路径，不是全局最优路径。

近期提出的等水平方法[27]是一种加工路径全局优化方法。与传统方法不同，其采用隐式加工路径表示方法(图 11)，并将各加工指标(如等残高，最大带宽)表达为能量函数，然后通过最小化这些能量函数获得全局最优加工路径，如全局优化残高分布[27]、光滑度[28]、长度[25]等。这一方法与场方法的结合具有很大的发展潜力。

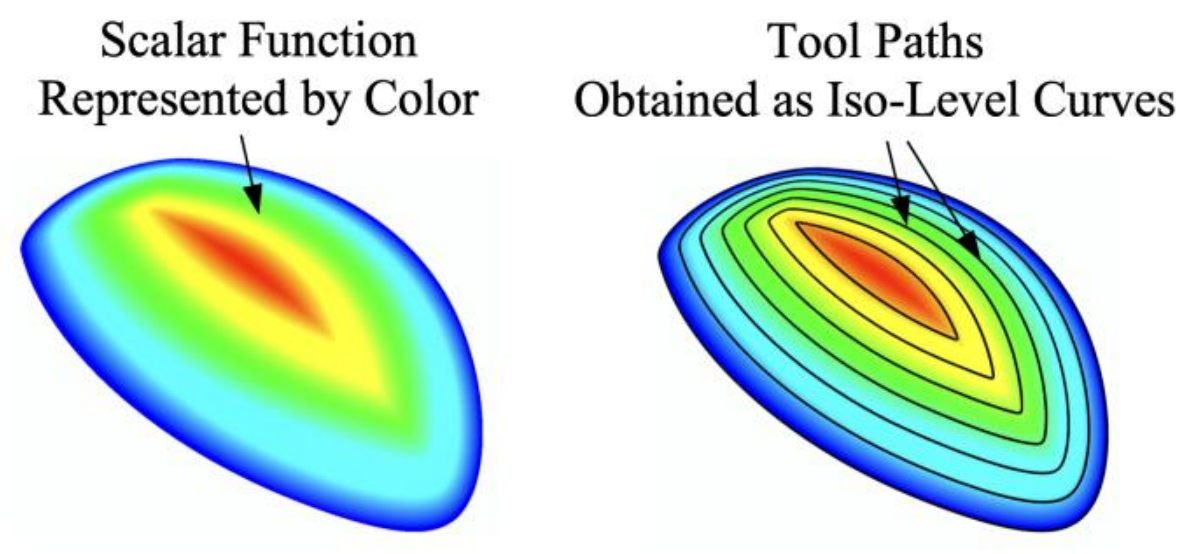

**图 11**　加工路径隐式表示示意[27]

除上述方法外，基于拓扑形式的空间填充曲线、摆线等路径生成方法[28]也得到了人们的广泛研究。该类方法在自由曲面抛光、恒定材料去除率(material removal rate，MRR)切削等方面有较好地应用。另外，人们还将启发式算法，如基因算法[29]、粒子群算法等，引入到加工路径全局优化中。这部分工作虽然仍处于初期，但具有很大地发展潜力。

## 2　刀轴方向优化

相比于三轴数控机床，五轴联动数控机床具有 2 个旋转轴(图 12)，从而可以加工三轴机床所不能加工的工件，同时能够以更高的效率和精度加工待加工曲面并提高表面质量。然而旋转轴的存在导致在五轴数控加工中出现工件坐标系和机床坐标系非线性耦合的问题[30]，使得轨迹规划比三轴加工时更为复杂。在五轴数控加工中除了需要规划刀位点的轨迹外，还需要对刀具姿态(刀轴方向)进行规划。

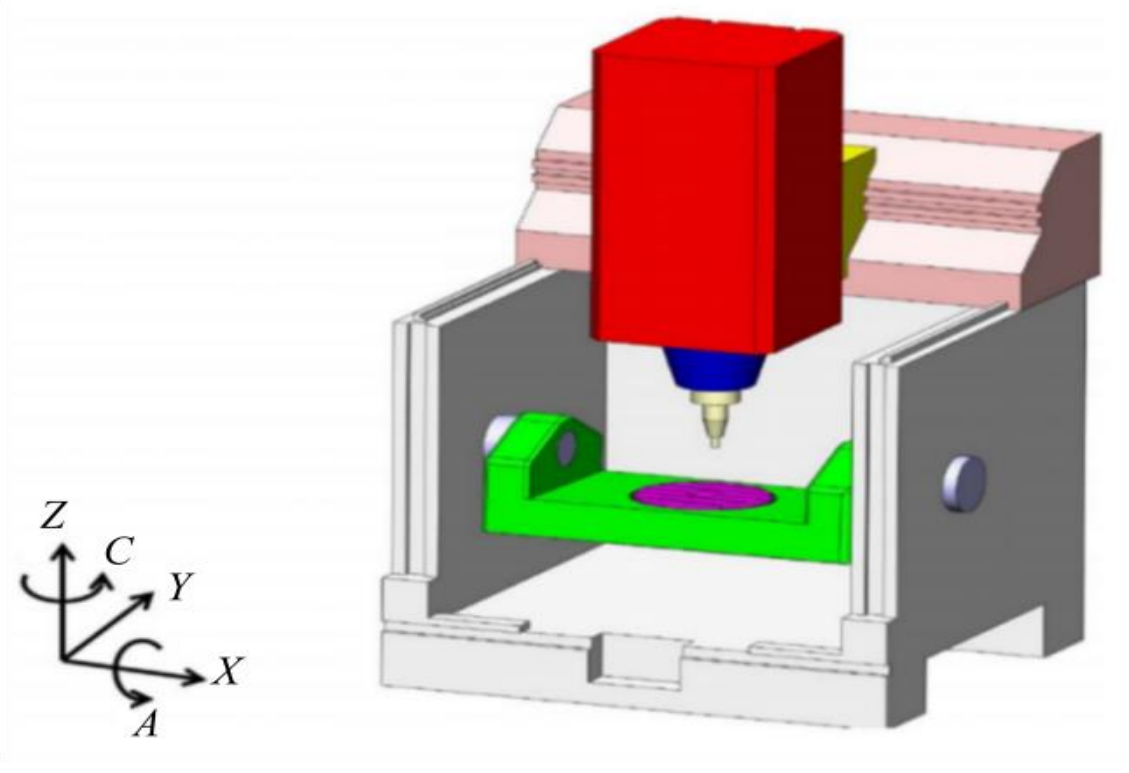

图 12　摇篮式五轴数控机床示意[30]

刀轴规划一般需要实现 2 个目标：避免干涉和刀轴方向的光顺。其中避免刀具干涉是基础，因为刀具干涉将会产生刀具和工件或者周围环境的碰撞，根据碰撞所产生的位置，大致可以分为全局干涉、局部干涉和尾部干涉(图 13)等情况(以下统称为干涉)。刀具干涉会损坏机床或工件，严重的会影响生产安全。无论哪种干涉情况，都会对刀具和工件造成不可修复的损害，因此必须保证高速加工中的刀具与工件之间是无干涉的。另一方面，在加工曲面曲率较大的地方时，可能出现刀轴方向突变的情况(图 14)，同样会影响加工效率和质量。因此需要对刀轴方向进行优化，保证整体刀轴方向是光顺的，得到更高地加工效率和加工质量。

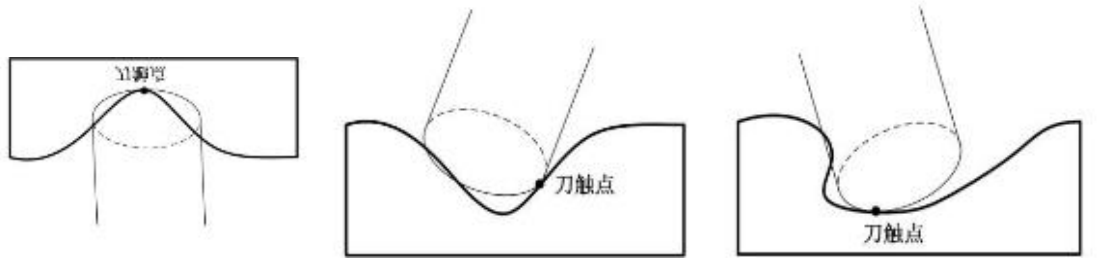

(a)　(b)　(c)

图 13　五轴加工中的三种干涉情况((a)局部干涉；(b)尾部干涉；(c)全局干涉)

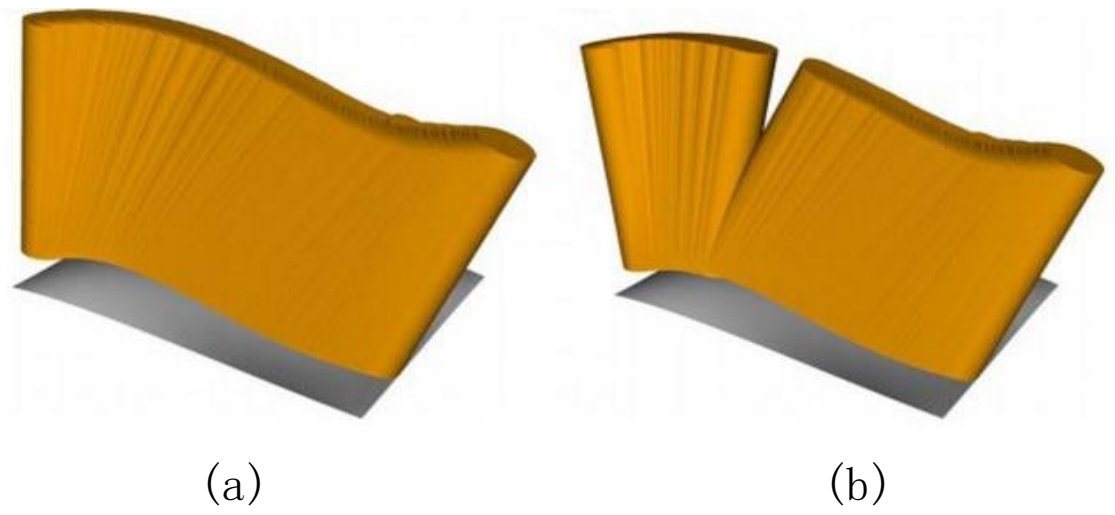

(a)　(b)

图 14　刀轴方向变化示意图[31]((a)光滑刀轴变化；(b)刀轴方向突变)

目前，刀轴优化的方法主要分为局部调整优化

和全局优化。其中局部调整方法的思路是在设计好的刀轴附近作微小的调整，使得刀轴方向具有更好的连续性。而全局优化的方法则是在刀轴方向的无干涉区域内进行优化，主要是基于构型空间(C空间或V空间)的方法。

## 2.1 刀轴方向的局部优化方法

### 2.1.1 考虑误差和刀轴方向的连续性

针对五轴折线段光滑过渡的工作主要考虑折线段连接处的光滑过渡处理。BI 等[32]在机床坐标系中用一条三次 Bézier 曲线光滑移动轴路径的拐角，再用一条 Bézier 曲线光滑转动轴路径的拐角(图15)。2 条曲线需要进行参数统一化处理，并能保证刀尖点和刀轴方向在工件坐标系中的误差。TULSYAN 和 ALTINTAS[33]提出了一个解耦的方法，在刀尖点和刀轴方向路径的拐角处分别插入五次样条，且在工件坐标系下考虑允许加工误差的约束，同时采用牛顿-辛普森迭代优化方法实现刀轴方向的光滑化和参数统一问题。YANGT 和 YUEN[34]在拐角处生成五次样条光滑路径，刀尖点路径在工件坐标系中表示，而刀轴方向路径在机床坐标系中表示，虽然实现了 2 条路径参数的统一和误差约束，但不能得到曲率极值点的解析解。HUANG 等[35]提出了一个解析的拐角光滑化算法，刀尖点路径在工件坐标系下进行规划，刀轴方向在机床坐标系下进行光滑化，根据加工时间将 2 条曲线进行参数化统一。

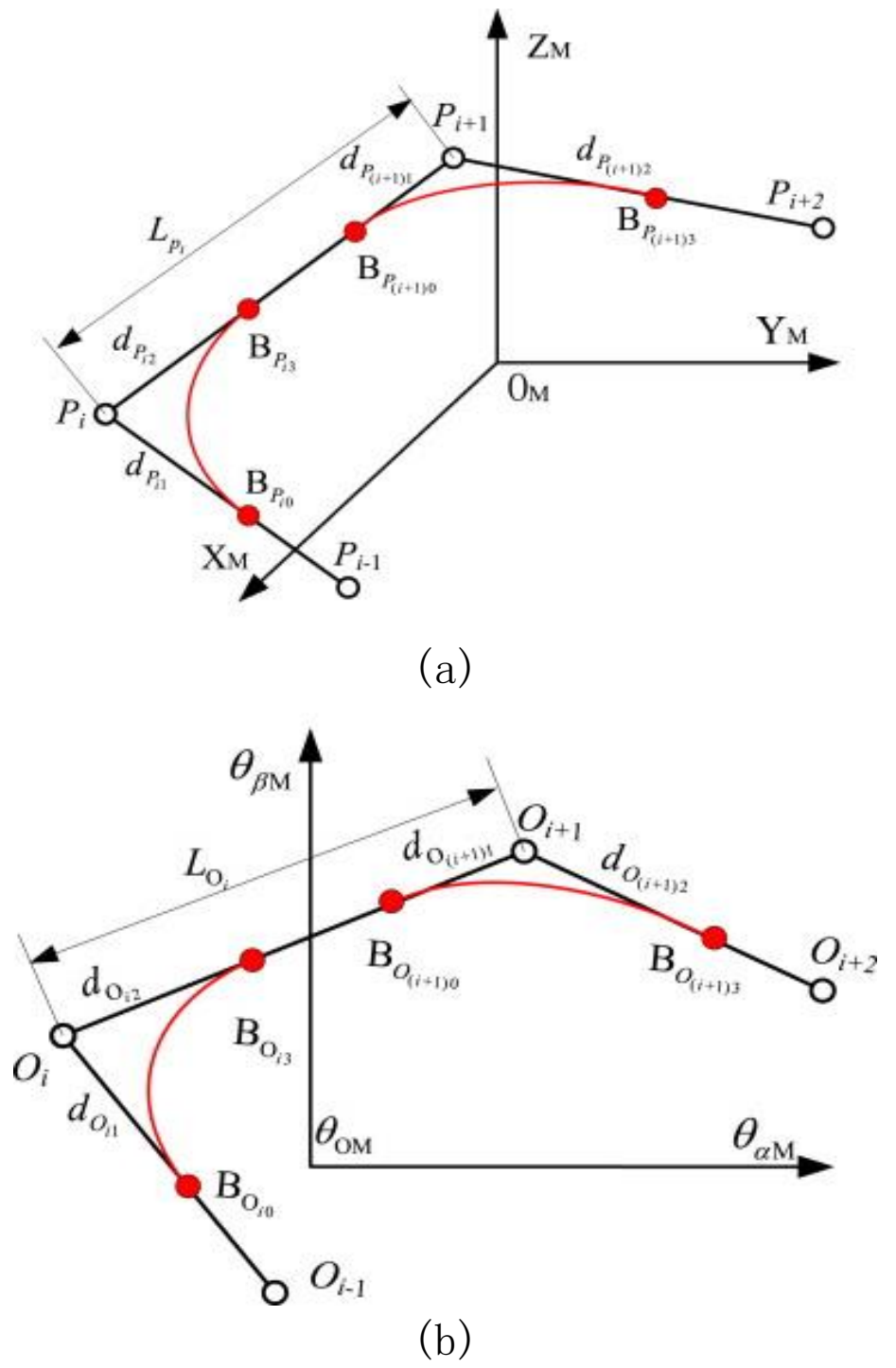

(a)

(b)

图 15 机床坐标系移动轴和转动轴拐角过渡曲线示意图[32]((a)线性移动轴拐角过渡曲线图; (b)转动轴拐角过渡曲线图)

### 2.1.2 刀轴的双样条表示和拐角光滑

此类方法用刀轴上 2 个点拟合得到的上下 2 条曲线来光滑化刀轴方向和刀尖点位置，以进行拐角处光滑处理。BEUDAERT 等[36]提出拐角处用 2 条三次 B 样条光滑过渡的方法，一条底部曲线定义刀尖点的位置，另一条顶部曲线定义刀轴上的第 2 个点(图 16)。用这 2 条曲线同时表示刀尖点位置和刀轴方向。并利用牛顿-辛普森迭代算法光滑刀轴方向的连接，得到较高的加工速度。SHI 等[37]在工件坐标系路径的拐角处插入一对满足刀尖点位置和刀轴方向误差的 PH 曲线，并通过优化控制点使刀轴方向的连续度更高。XU 和 SUN[38]在拐角处获得更小的曲率和更高的加工速率，但在计算曲率极值点时计算耗时较大。

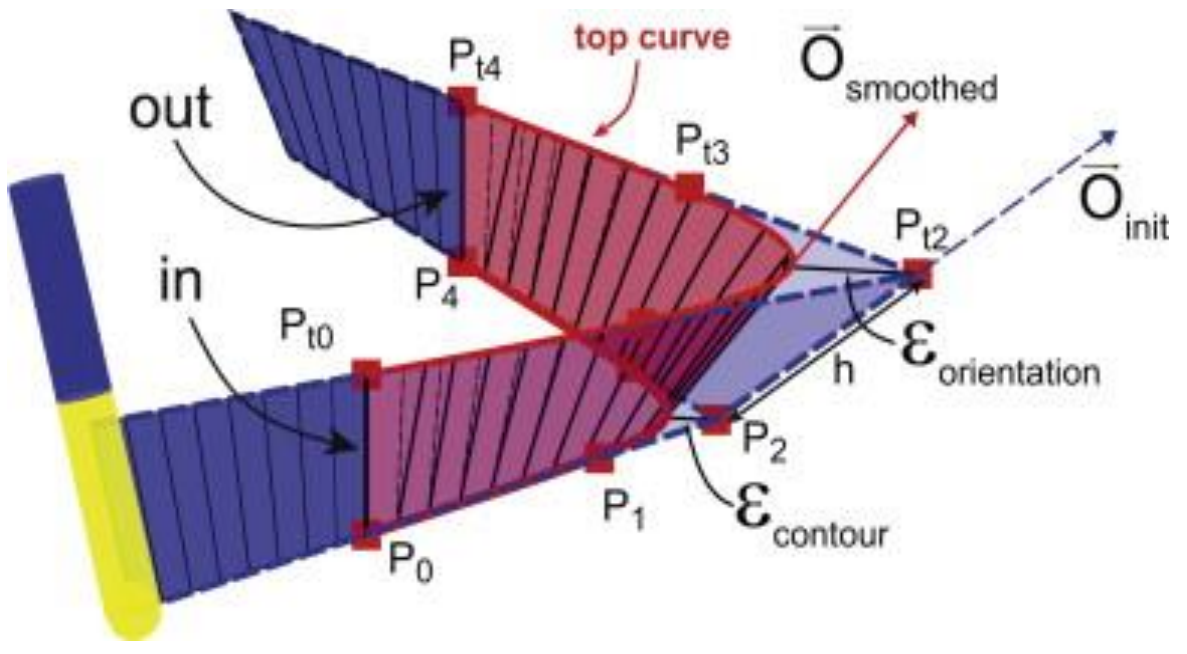

图 16 五轴加工拐角过渡及误差示意图[36]

## 2.2 基于 C 空间的全局刀轴方向优化方法

以下以摇篮式机床为例介绍基于 C 空间的刀轴优化方法。对于五轴数控加工，在局部坐标系中刀轴方向完全由俯仰角 $\theta$ 和倾斜角 $\phi$ 确定(图 17)。将三维空间中单位球的上半球映射到由俯仰角 $\theta$ 和倾斜角 $\phi$ 张成的二维平面，则平面上每个点都代表了空间中的一个刀轴方向。

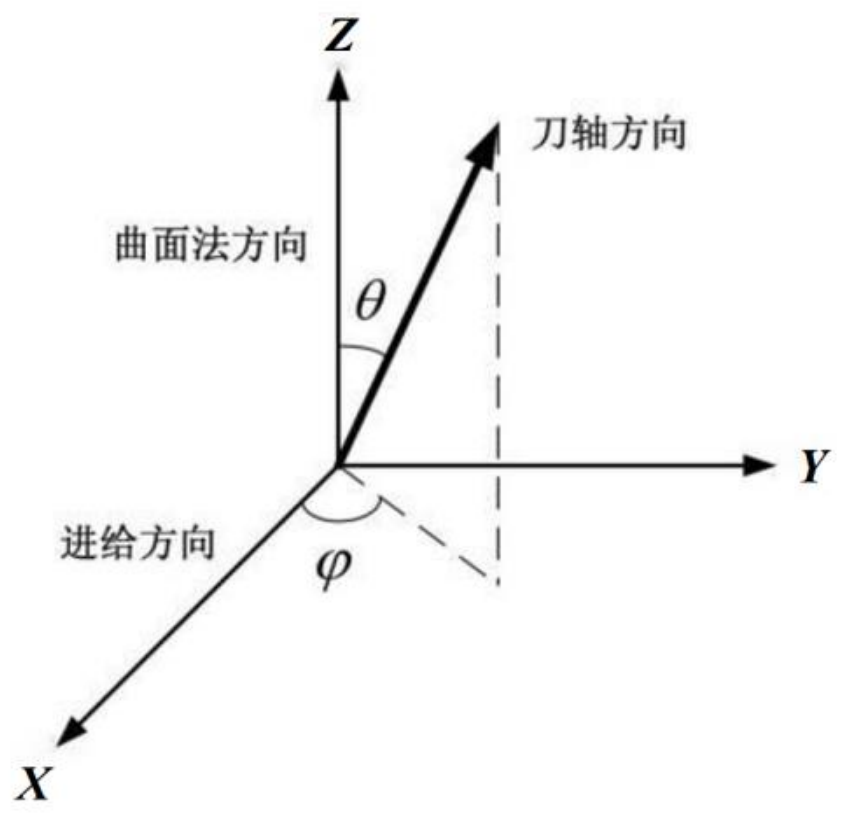

图 17　局部坐标系与刀轴方向

在每个刀位点处,所有可行(无干涉)的刀轴方向组成的集合称为该点处的构型空间，或称 C 空间。为了在 C 空间内规划刀具姿态，需要先算出每个刀位点处的 C 空间。但是，计算曲面上某点处的精确 C 空间是非常复杂和困难的。对于球头刀，MORISHIGE 等[39]将加工工件视为多面体，计算了球头刀与多面体的碰撞边界，进而计算出二维 C 空间(图 18)，最后将不可行的刀轴方向移动到 C 空间的边界，从而生成无碰撞的刀轴方向。关于平底刀和环形刀的工作此处不再赘述。

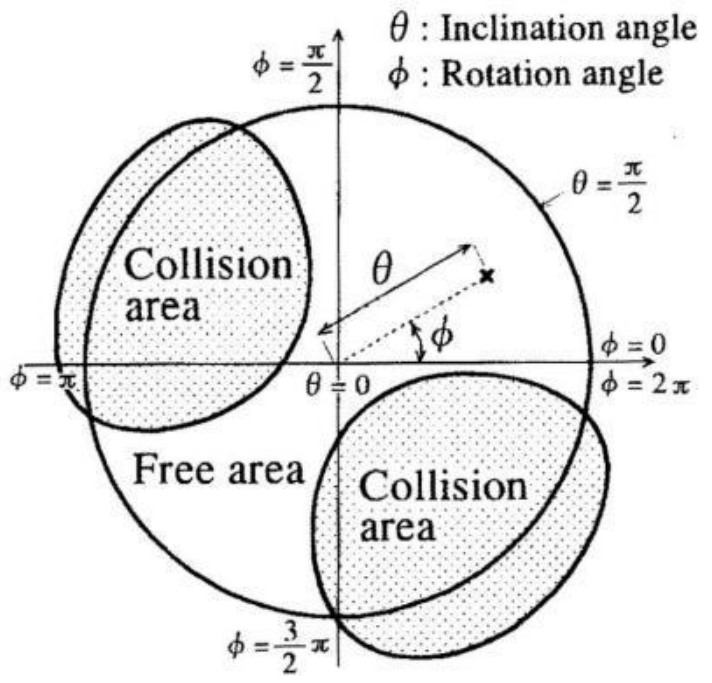

图 18　2 维 C 空间示意图[39]

求得每个刀位点的 C 空间后，刀位点序列的 C 空间就构成了一个三维的可行区域，刀轴优化问题就是在这个三维可行区域中选取一条连续的路径，使得刀轴光顺，利于高速高精加工。这方面有着丰富的研究工作。JUN 等[40]提出一种搜索算法来计算 C 空间，并通过极小化刀轴突变来光顺刀轴方向，能够得到较好的结果。WANG 和 TANG[41]在限制刀轴变化角速度的情况下，通过构建可视图(visibility map，V-map)来避免干涉问题，保证了刀轴角速度变化的一致性。LU 等[42]考虑了抬刀高度,提出一种基于三维 C 空间的刀轴规划方法，可同时考虑残高和干涉情况。LIN 等[43]在 C 空间中整体调整刀轴方向，以避免刀轴方向的奇异位置。

在只考虑刀轴轨迹的光顺性方面，HO 等[44]利用四元数插值优化刀轴轨迹。CASTAGNETTI 等[45]给出了可行方向域(domain of admissible orientation, DAO)的概念，将机床各轴的速度、加速度的限制融入刀具姿态的优化，进一步改善了刀具的运动性能。YE 等[46]针对材料去除率和刀轴的二阶运动学的连续性，将刀轴优化问题转化为二阶运动学约束下的时间最优插补的优化问题进行求解，得到了兼顾加工效率和加工质量的刀轴优化方法。文献[31]将可行域进行离散后，利用图论中的算法，生成比较光顺的刀轴轨迹。文献[47-48]为了得到更为光顺的加工路径，引入了差分图的概念，得到二阶光滑的刀轴路径。这一方法适用于任意刀具，区别在于 C 空间的计算方法不同。该方法在求得每个刀位点的 C 空间后，将其离散化，得到每个刀位点处离散的可行刀轴方向。然后以可行刀轴方向为顶点，将相邻 2 个可行域内的顶点进行连边得到一个分层有向图，每条边以前后 2 个顶点的角度变化值赋予权重。利用 Dijkstra 算法，能够找到图中的最短路径，并找到最短路径上各顶点对应的刀轴方向，即为最终刀轴轨迹。这一方法的优势在于得到全局优化的刀轴路径的同时，通过修改 C 空间的结构，较容易避开机床的奇异位置。以分层图的每个边为顶点，以前后 2 个边的角度变化的差(可以看成刀轴方向的二阶微分)为权，构建一个差分图。对其进行最短路径求解，就可以得到最终兼顾角度的一阶和二阶光滑度的刀轴优化方向。ZHANG 等[49]应用强化学习方法给出了一个刀轴方向的全局优化方法。

## 3　G代码加工

在数控加工中，为了加工复杂曲面，通常使用小线段去逼近刀具的运动轨迹，因此小线段(G01)插补是数控加工中的核心问题之一。如果机床严格按照小线段所表示的路径进行加工，为了满足加工误差，在线段拐角处的允许速度非常低，引起频繁的加减速，降低加工的质量和效率。因此，目前小线段插补的方式主要分为以下几类：

(1)　小线段拐角过渡插补，即在 2 个小线段处插入一段光滑曲线，达到光滑化刀具路径的效果；

(2)　小线段拟合成光滑曲线，再对光滑曲线进行插补；

(3) 基于 FIR 滤波的过渡插补算法。

### 3.1 小线段拐角过渡方法

小线段拐角过渡方法分为 2 种。

第一种方式是根据加工误差约束在拐角处插入一条曲线(图 19)，然后整体进行小线段和拐角处曲线的速度规划。HAN 等[50]在拐角处插入三次 B 样条。SUN 和 ALTINTAS[51]在拐角处插入 G3 连续的双 Bézier 样条曲线，满足加加速度有界和误差约束。ZHANG 等[52]使用 G4 连续的 B 样条曲线在拐角处过渡，考虑了切向加速度以及加工误差约束。以上方法所采用的这种方式虽然可以满足加工误差，也可以提高拐角处的加工速度，但是由于需要进行曲线的速度规划，不能满足实时加工的要求。

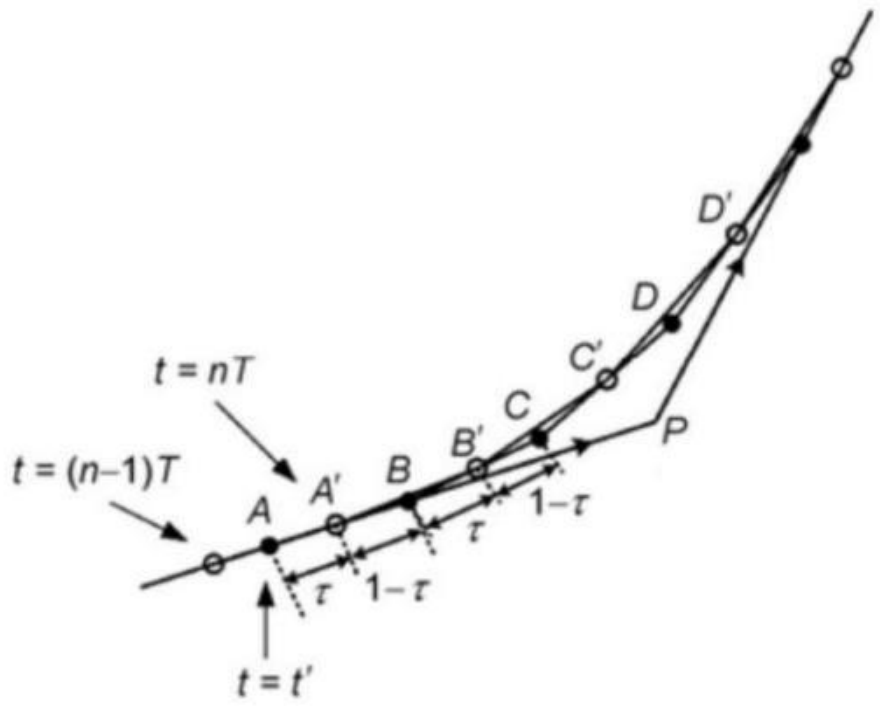

图 19　折线段拐角过渡示意图[52]

第二种方式是根据加工误差和拐角处局部满足 Bang-Bang 最优控制为目标，确定一条以时间为参数的拐角过渡曲线。该曲线自动满足加工误差、机床动力性能约束且局部加工时间最短为优化目标的数学模型。该方法另外一个最大的优势是自动生成速度规划曲线，满足第一种方式所不具备的高效计算，达到数控系统的实时加工要求。这种方式又分为加速度有界的速度规划、加加速度有界、jounce 有界，三角函数速度规划方法等。DUAN 和 OKWUDIRE[53]在拐角处插入以时间为参数的 NURBS 曲线。LIN 等[54]提出一种基于加加速度有界的拐角过渡方法，该方法同时考虑了机床各轴运动学性能和刀尖点的几何误差，不过未考虑转动轴的几何误差。LI 等[55]根据 Sigmoid 函数的特点，使用其作为拐角过渡曲线，满足切向加速度以及几何误差约束。ZHANG 等[56]在五轴数控加工时，在局部拐角处根据 Bang-Bang 最优控制目标，生成满足加速度有界的拐角过渡曲线，根据整体前瞻算法，进行全局速度优化。该算法适用于五轴数控加工，并将加工路径表示在五维空间中，通过将工件坐标系中拐角误差转化到机床坐标系，实现了误差与机床各轴运动学性能的约束，且计算效率满足在线加工的要求。

### 3.2 小线段拟合与插补方法

在三轴数控加工时，常用 B 样条拟合刀具路径，B 样条拟合一般采用最小二乘拟合方法和前进迭代逼近算法，最小二乘方法只考虑数据点间误差，未考虑中间过程的误差。这些方法得到的 B 样条不以弧长为参数，CHEN 和 KHAN[57]使用最小二乘方法将 B 样条转化为以弧长为参数的曲线，但是折线与曲线间的误差并没有严格限制。在五轴数控加工时，刀位点和刀轴方向矢量需要同时进行拟合，相对于三轴数控加工，曲线拟合复杂度提高。目前已有的方法主要采用 B 样条曲线、球面样条曲线[58]和四元数曲线这 3 种方法。MIN 等[59]在六维空间对五轴的加工数据进行 B 样条拟合，满足刀心点和刀位点以及刀轴方向的误差约束。

然而，以上方法生成的曲线刀具路径，在数控系统中需要进行速度规划，转化为以时间为参数的插补点序列，该计算过程相对耗时，对于高次参数曲线较难达到实时加工的要求。为此，ZHANG 等[60]将 G01 代码表示的折线段利用二次样条曲线进行拟合，然后通过分轴加速度有界的时间最优插补方式进行插补计算[61]，并给出了速度函数的显示解，最终得到了误差和机床运动学约束下的高效加工方式，优化了加工路径和加工效率。

### 3.3 基于 FIR 滤波的小线段过渡插补算法

有限脉冲响应 (finite impulse response，FIR)滤波器插补算法是实时的拐角过渡算法，其特点为不需要先确定拐角路径、然后规划拐角速度 2 个步骤，只需根据拐角误差约束自动规划拐角路径和拐角速度。首先，规划点到点的速度曲线时，匀速运动曲线经过 2 个线性滤波器卷积变成光滑的 7 段式速度曲线，可根据线段长度、加速度和加加速度约束等求得 2 个滤波器的延迟时间，从而确定具体的速度曲线，如图 20 所示。然后，规划连续路径，为了提高运动光滑度和加工效率，将相邻 2 个线段的速度曲线复合一部分，即前一段线段的运动还未完成时，提前开始下一段运动，这会导致拐角误差，根据拐角误差约束确定速度曲线复合段的时间，拐角速度和拐角误差都可由相邻 2 段运动的复合计算得到，从而完成连续线性路径的插补。

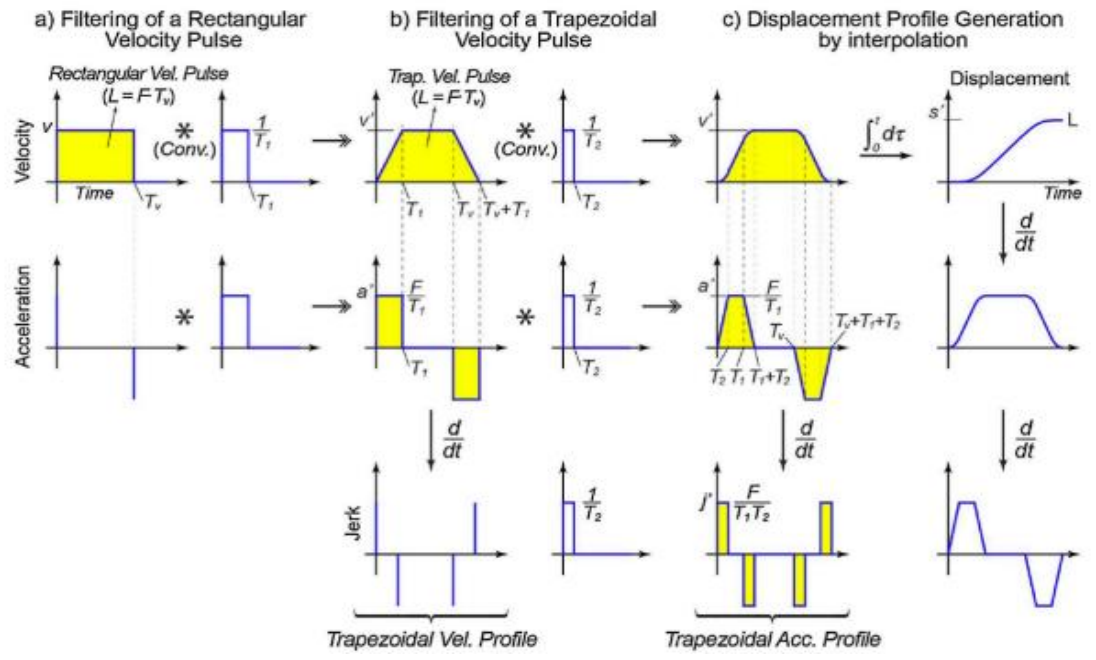

图 20　基于 FIR 滤波器生成的点到点的速度曲线[57]

BIAGIOTTI 和 MELCHIORRI[62]提出用 FIR 滤波器规划加加速度约束的点到点的速度。文献[63-64]提出用 FIR 滤波器插补三轴连续线性路径的实时算法。文献[65-66]提出用 FIR 滤波器插补五轴连续线性路径的实时算法。在五轴 FIR 滤波器插补算法中，滤波器分别用于刀位点和球形坐标表示的刀轴方向，通过统一用于刀位点和刀轴方向的滤波器的时间常数，实现各刀轴的运动统一。

除了上述几类算法，还有其他相关的算法。如通过结合拐角过渡与拟合插补 2 种方式，对小线段进行分段，分为长直线段和连续的微小折线段，在连续的微小折线段处进行曲线拟合插补，在长直线段和微小折线段处进行拐角过渡处理。此外，针对特别微小的折线段，拐角过渡并不合适，YE 等[67]针对这一问题，根据其几何特征进行分段和速度标记，然后通过对分段的折线段进行加速度有界约束下的前瞻处理使得得到的速度可达，从而使得插补可以跨越多个微小折线段。

# 4　参数曲线的进给速度规划方法

进给速度规划是加工轨迹插补中最重要的过程之一，其目的在于在几何和运动学条件约束下生成平滑的刀具运动进给率曲线，对加工过程特别是加工效率和质量方面具有重要影响。

在实际加工轨迹插补过程中，如何在满足限制条件下，尽可能得到具有更高加工效率的进给率值，成为进给率规划中重要的研究内容。主要的限制可以分为三类：

(1)　弓高误差限制：相邻插补点之间的直线与原曲线之间的偏差；

(2)　轮廓误差限制：实际加工轨迹与理想加工轨迹之间的偏离；

(3)　加速度和加加速度限制：反映了驱动电机所能提供的最大扭矩以及提供扭矩变化的能力。

常用的进给速度规划方法可分为加减速模型方法、优化方法和光滑进给一体化方法 3 类。加减速模型方法通常给定一个固定的函数类型，用来实现 2 个不同速度阶段之间的平滑过渡，常用的函数类型包括梯形函数、S 型函数、Sigmoid 函数和三角函数等。对于给定的刀具轨迹，首先会选出一些间断点和曲率极值点作为关键点，进而根据几何与运动学的限制给出关键点处进给速度的上限值，并最终利用给定的函数类型，实现相邻关键点之间的速度变化衔接。优化方法一般未给定固定的函数类型，通常使用一条或多条自由曲线来表示速度进给曲线。对于给定的刀具轨迹，通过采样一些刀具轨迹点，以进给速度最大为目标，以几何误差和运动学约束为限制条件构造优化问题，并最终求解优化问题，得到刀具轨迹对应的进给速度曲线。光滑进给一体化方法则将刀具轨迹的光滑化处理和进给速度规划合并为一步处理。经典方法包括有限滤波(FIR)方法、运动交叠方法和时间样条方法等。

## 4.1　加减速模型方法

由于加速度的剧烈变化会导致加加速度过大，从而引起机床的不良振动，目前的加减速模型通常考虑让加加速度小于某个给定阈值，这个阈值是由机床特性决定的。ERKORKMAZ 等[68]首先提出使用 $S$ 形函数来生成加加速度受限的速度曲线，如图 21 所示，在这个 $S$ 形曲线中，加加速度为阶梯状，且取值为最大值、零和最小值这三者之一。该 $S$ 形曲线由 7 段组成，加速阶段 3 段，匀速阶段 1 段，减速阶段 3 段组成。该 $S$ 形曲线能实现 2 个不同速度之间的平滑过渡。在 JAHANPOUR 和 ALIZADEH[69]提出的方法中，一个 $C^2$ 连续的 5 次 Bézier 函数被用来生成加速度和加加速度均连续的速度曲线。JIA 等[70]提出了速度敏感区域的概念，在速度敏感区保持匀速，在过渡区平滑过渡不同速度敏感区的进给速度，减少了进给速度的变化频率，提升了加工质量。

为了生成更光滑的进给速度曲线，一些加加速度连续、加加速度的导数(简称震动)受限的模型被提出。FAN 等[71]提出了由 15 段构成的进给速度模型，其中加速阶段 7 段，匀速阶段 1 段，减速阶段 7 段。该速度模型中，震动为阶梯状，且取值为最大值、零和最小值这三者之一。ZHANG 和 DU[72]对

上一模型进行了简化，删掉了加(减)速阶段中加加速度恒定的 2 段，使得相关时间参数的计算过程得以简化，但是总体插补时间有小幅度增加。WANG 等[73]将 7 段的 $S$ 形函数简化为 5 段，即不存在加速度恒定阶段，并使用三角函数代替阶梯函数，实现了加加速度的连续。由于该模型的分段数目大大减少，相关时间参数的计算过程被简化。但是由于三角函数中，加加速度仅在少数几个点达到极大值或极小值，总体插补时间也会有小幅度增加。

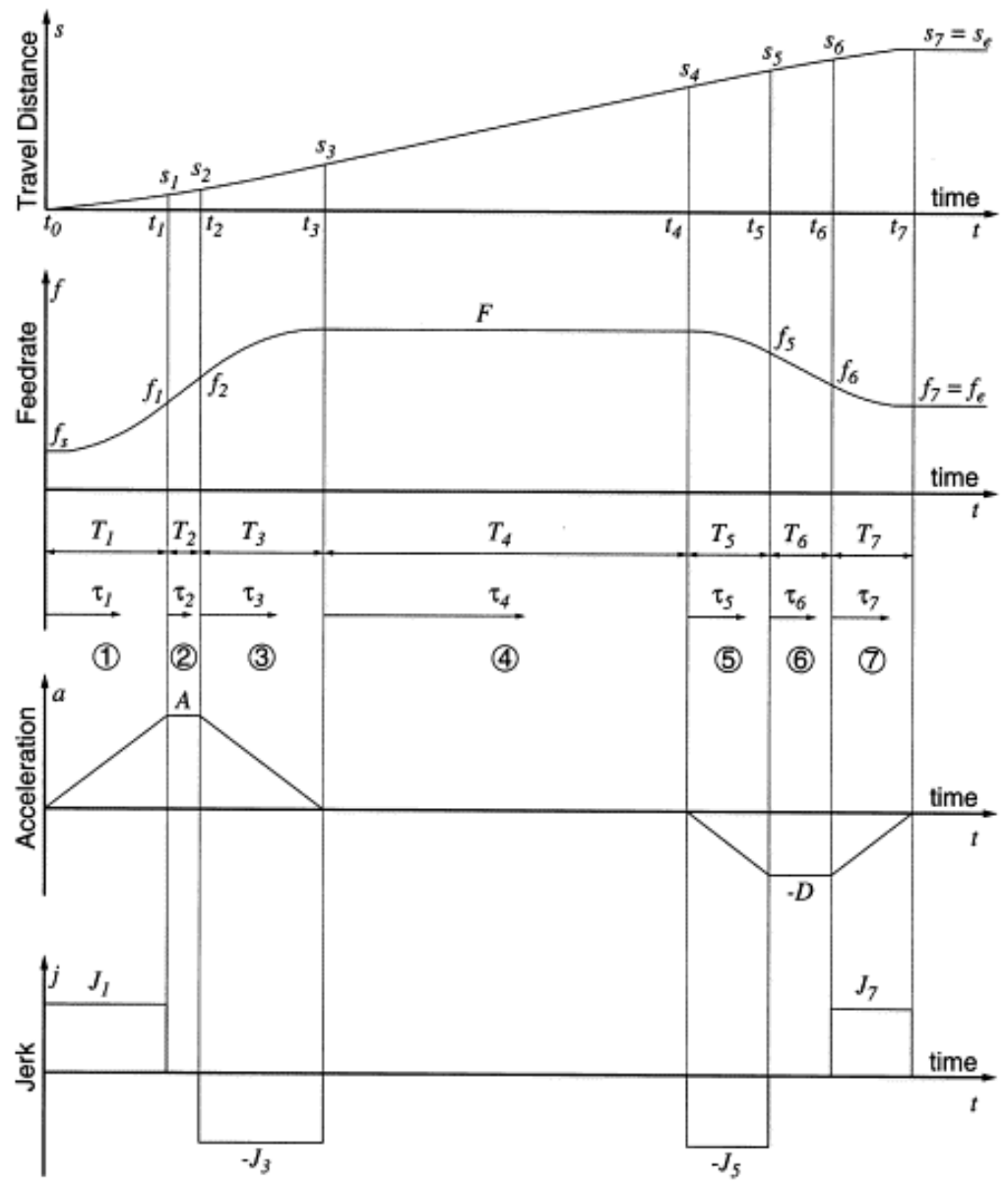

图 21　S 形函数的运动学轮廓[68]

## 4.2　优化方法

SUN 等[74]提出了一种考虑几何、工艺和驱动约束的五轴加工进给速度调度方法。首先利用限制弦误差、角速度和轴速度构建初始进给速度曲线。然后，对进给速度敏感区域进行比例调整，以同时减小角加速度、线加速度、轴加速度和急动等约束的大小(图 22)。由于涉及迭代调整，该方案无法实时实现。后 SUN 等[75]又提出了一种对相邻 B 样条进给速度曲线具有 $C^2$ 连续性保证的分段线性规划方案，以实时解决进给速度调度优化模型。然而，在上述方法中，五轴刀具轨迹曲线是用归一化弧长参数进行参数化的，这是一个相对较强的条件，需要提前对刀具轨迹弧长参数化。此外，针对双机器人对雕刻曲面零件的同步加工，SUN 等[76]建立了双刀具路径上对应点之间的同步映射关系，求出保证两刀同步运动的约束条件。然后，建立了具有加加速度和同步约束的时间最优同步进给调度模型，该模型可以充分利用机器人的运动学性能，并充分保持刀具的切削性能。为了进一步加快非线性模型的求解过程，构造约束相对于进给速度平方的线性关系，从而使用线性规划算法来实现进给速度曲线的精确表达。

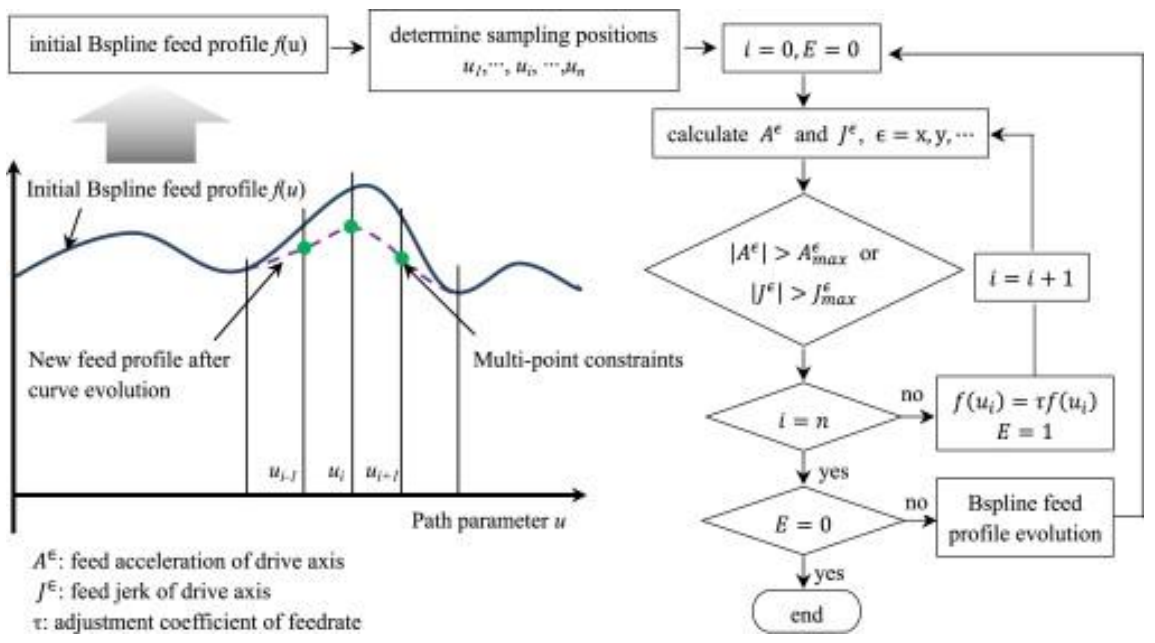

图 22　比例调整的进给速度优化方法[75]

LIU 等[77]通过使用原始连续目标函数的离散格式，推导出了用于优化进给速度的线性目标函数。然后，将预设的多约束近似为线性目标函数中决策变量的线性约束条件。进给速度优化的线性模型可以通过成熟的线性规划算法进行求解，最优解被拟合到光滑的样条曲线上，作为最终的进给速度曲线。为了保证弗勒内-塞雷(Frenet-Serret)框架下各种约束均被满足，YE 等[78]首先将参数曲线形式的刀具轨迹弧长参数化，并分析了在 Frenet-Serret 框架下速度、加速度和加加速度的约束。然后引入时间最优方法在时域生成采样点，并将时间最短为优化目标的优化模型转换为以速度的平方和最大为优化目标的优化模型，采用线性方程求解该模型的次优解，并采用三次 B 样条拟合采样点处的速度，最终得到进给速度曲线。

YUAN 等[79]提出弓高误差、进给速度和切向加速度界下的时间最优插补算法和弓高误差和切向加速度界下的贪心插补算法。其核心思想是将弓高误差限制在向心加速度限制上，从而得到速度极限曲线，称为弦误差速度极限曲线。然后，速度规划是在弦误差速度极限曲线下，找到由加速度或加加速度限制的适当速度曲线，如图 23 所示。

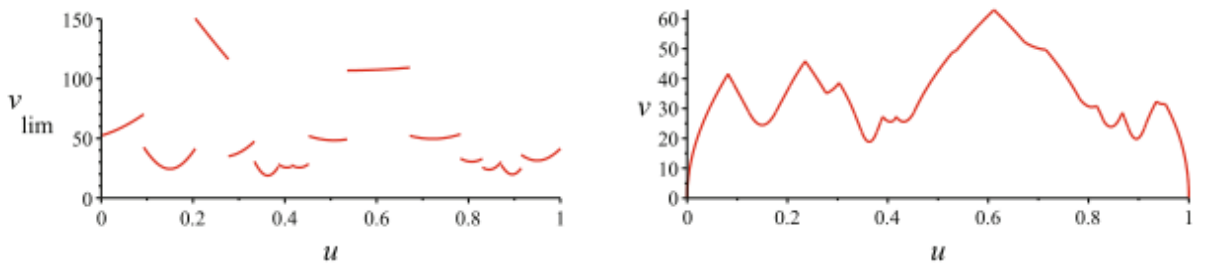

(a)　　(b)

图 23　速度极限及其构成的速度函数[79]((a)弦误差速度极限曲线；(b)进给速度曲线)

### 4.3 光滑进给一体化方法

LIU 等[80]提出了一种在工业五轴机床中连续插补 G01 指令的方法。采用有限脉冲响应(finite impulse response，FIR)滤波器来生成线性段的时间最优和加加速度限制的轨迹,其中平移运动和旋转运动分别在 WCS 和 MCS 中进行过滤。并建立了一个缩放原则使平移和旋转运动相对于 FIR 滤波器的持续时间同步，同时最小化循环时间。与现有的两步几何方法不同，G01 指令通过一步直接重叠相邻的运动学曲线来连续插值。基于所提出的容许面积概念，通过调制相邻 2 段的重叠时间，将拐角误差限制在给定范围内。为了确保轴向和切向运动学约束均得到很好的保持，提出了一种新的预测和校正方案来预测平移运动的最大可达切向运动学约束，同时校正最大可达到的旋转运动的轴向运动学约束。TAJIMA 和 SENCER[81]提出了一种新颖的在线奇点避免的 5 轴机床实时轨迹生成技术。其中平移和旋转刀具运动由进给脉冲信号表示，并经过低通滤波以产生具有受控频谱的同步平滑运动。过渡误差通过调节进给速度来控制，降低进给率可减少过渡错误。因为路径混合控制和进给率规划是共同实现的，这为进给速度规划打开了新的潜力。该算法，通过直接卷积 2 个串联的 FIR 滤波器链来保证进给速度曲线的全局切向平滑度。为了控制过渡误差，分析了滤波诱导误差的原理，该误差可以通过进给参数和曲率半径来计算，并在轴向运动学约束之外增加了最大误差作为额外的约束。

ZHANG 等[82]使用一步过渡轨迹生成策略，其以时间为参数的二次曲线，且在一定意义下是时间最优的。此外，该算法还考虑了合理的轴向加速度限制，并提出了一种新的拐角过渡模型，以充分利用拐角误差，过渡路径是由最大轴向加速度和最大误差容限确定的时间参数化二次曲线。TAJIMA 和 SENCER[83]提出了考虑加加速度限制的对称运动学转角平滑方法，该方法直接规划了在转角处的加加速度限制的进给率曲线实现局部转角光滑。并基于局部光滑方法针对连续微小线段提出了全局的光滑方法，通过对局部光滑进给率延长或合并，在满足角度误差的前提下实现全局的光滑。ZHANG 等[84]使用对称的加加速度限制模型进行每段直线轨迹的进给率规划，然后对转角过程运动交叠的几种情况进行分析，将轨迹光滑问题转化为运动时间优化问题，在满足转角误差、各轴加速度和加加速度的限制条件下，计算每段直线上的切向加加速度，随后将这种方法推广到五轴的情形上。为了充分利用各轴的运动特性，WANG 等[85]提出了一种局部非对称运动角平滑方法，其中局部转角光滑的进给速度是非对称的，这意味着减速和加速过程的位移可以按照相邻直线段的位移大小进行调节。

随着人工智能技术的蓬勃发展，相关理论与方法也陆续被应用到进给速度规划中。如粒子群优化方法[86]和遗传算法[87]。人工智能方法在多目标和多约束优化中具有天然的优势。即：在目标函数中增加相应的项，就可以起到对多个不同目标的约束；在模型中增加相应的条件，就可以对高阶、非线性的约束起到良好的限制效果。人工智能方法在数控加工中的实时应用，将会是未来发展的一个重要方向。

## 5 新型加工优化

在数控机床的快速发展过程中，人们在实际加工中对各个流程的精度与效率提出了更高的要求，因而涌现出了众多的新型加工优化算法。下面本文将从 CAD，CNC 的 22 个方面介绍目前最新的一些数控加工方法。

### 5.1 CAD：满足加工约束的新型自由分片逼近

七十年代，计算机开始应用于汽车、舰船和航空、航天工业，由此催生了计算机辅助设计(CAD)这一领域。但是，早期的奠基者多是工程师和计算数学家，其设计的 CAD 理论和方法只考虑了拓扑简单的局部曲面设计，未考虑复杂拓扑曲面的整体设计，以及相应的加工特性约束，从而导致在后续的流程中，加工设计好的曲面，需要手动进行多项安全测试(如刀具防碰撞检测等)，从而间接降低加工效率。为此，有学者提出，在实际加工中，给定曲面进入到 CAM 环节之前，对其进行一系列的加工可行性分析。在满足误差约束的前提下，对曲面进行分片逼近，使得输出的结果可以直接用来加工。

在规划刀具路径之前，需要先配置一个刀具加工方向，此时待加工曲面上任意一点的外法向量与刀具加工方向的内积必须大于 0。所以在加工一个完整的复杂曲面时，由于数控机床加工空间受自由度的限制，几乎不可能仅靠一个刀具方向

来加工完成。且需要对这个复杂曲面进行分片处理，然后在每一片上进行路径规划并加工，如图 24 所示。

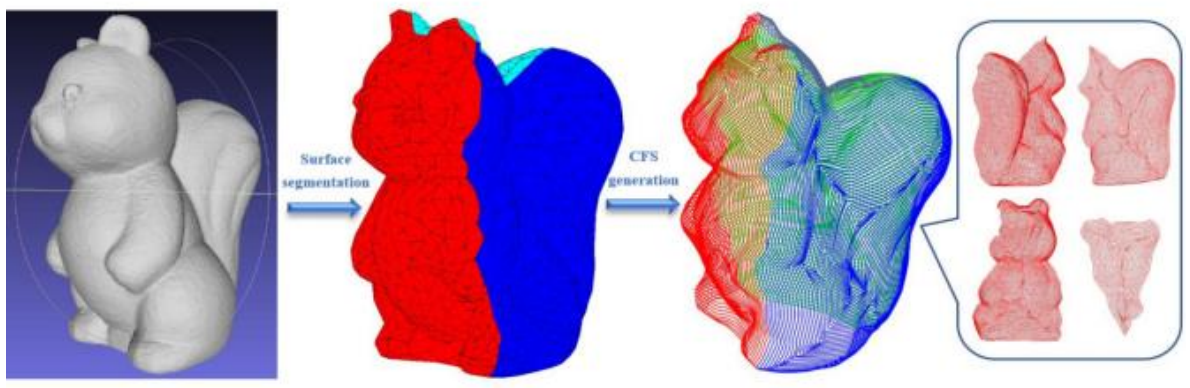

图 24　数控加工流程示意

通常情况下，上述如规划全局加工流程、保证刀具可达性等任务是交给富有经验的数控工程师来手动完成。然而，随着 3D 建模能力的增加以及由此产生的待加工曲面的复杂性使得手动任务愈加困难。因此，许多计算机图形学文献研究了如何将曲面或铣削体自动分解为适合数控加工的若干曲面片。

而基于高度域(Height Field)的曲面分片方法无论是在算法效果，还是发展潜力方面均最具代表性。此类方法中，输入的自由曲面的表达形式一般为三角网格模型 $M=(V, F)$。$V$ 表示三角网格的顶点集合以及各自的对应坐标。而网格内的三角形可表示其顶点索引 $F$。目标是将三角网格模型分割成少量的曲面片，以便每个曲面片都适合用于数控加工，且满足约束

$$E_i \subseteq \mathrm{F}, \bigcup_i E_i = \mathrm{F}, E_i \cap E_j = \emptyset,$$

其中，$E_i$为给定曲面分割出来的曲面片，如果满足以下所提的 5 个加工约束，则称之为高度域。

(1) 刀轴方向约束：对于每一个曲面片$E_i$，均要给定一个刀轴方向$\overrightarrow{d_i}$使整个$E_i$均可被加工到，也就是说，$E_i$上任意一点的外法向量与$\overrightarrow{d_i}$的内积均大于 0。

(2) 体积约束：大多数数控加工设备对毛胚都有体积限制，即长、宽、高不能不超过机床容许最大值[88]。而这个约束同样适用于各曲面片。

(3) 可达性(Accessibility)约束：在数控加工中，需要无碰撞地加工曲面。在待加工曲面上任意一点都有一个方向集合，且包含所有在该点加工时刀具与曲面无碰撞的方向。这个集合被称为可达性锥(Accessibility Cone)[89]。图 25 显示了二维平面下的锥体，而在三维空间该集合往往会形成一个圆锥体。一个三角形上所有点的可达性锥的交构成了这个三角形的可达性锥。一个曲面片$E_i$满足可达性约束，即为$E_i$上所有点的可达性锥的交不为空集。

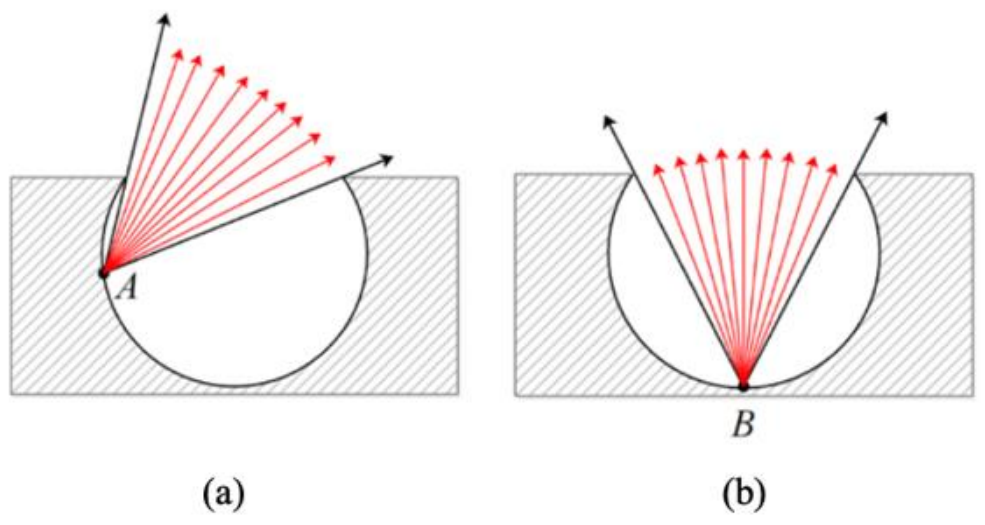

图 25　2 个点处的可达性锥(以二维平面为例)：用黑色实线表示的方向为可达性锥的边界((a)点 $A$ 处的可达性锥；(b)点 $B$ 处的可达性锥)

(4) 曲面片数量约束：在加工一个完整的复杂曲面时，由于数控机床加工空间自由度的限制，几乎不可能仅靠一个刀具方向完成加工。因此需要对这个复杂曲面进行分片处理，然后在每一片上进行路径规划并加工。每加工一个分块后的曲面片，就需要对刀具进行一次装卡定位。这样的操作费时费力，所以总希望装卡定位操作的次数尽量少，也就是说，分块的曲面片数量要尽量少。因此，应该将待加工自由表曲面分割成少量的曲面片，即为曲面片数量约束。

(5) 边界约束：在曲面分块后，需要在各个曲面片上生成刀具加工路径。而现有的很多刀具路径生成算法都会基于曲面边界做偏置生成路径轨迹。因此，曲面片边界的光滑性影响着刀具路径的质量。拥有太多拐点的边界可能会导致这个曲面片上生成的路径轨迹也充满了拐点，从而降低了进给速度，影响加工效率。因此，每个曲面片的边界必须足够光滑，这就是边界约束。

在给出相关约束后，高度域方法将曲面分片问题转换成一个离散的打标签(即对于每一个三角形，给定一个合理的方向标签)问题，并根据相关约束定义能量函数。之后，此类方法根据三角网格的拓扑连接关系建立加权连通图，如图 26 所示，并利用求解连通图的最小割来最小化能量函数，进而完成曲面的分片，得到各高度域。

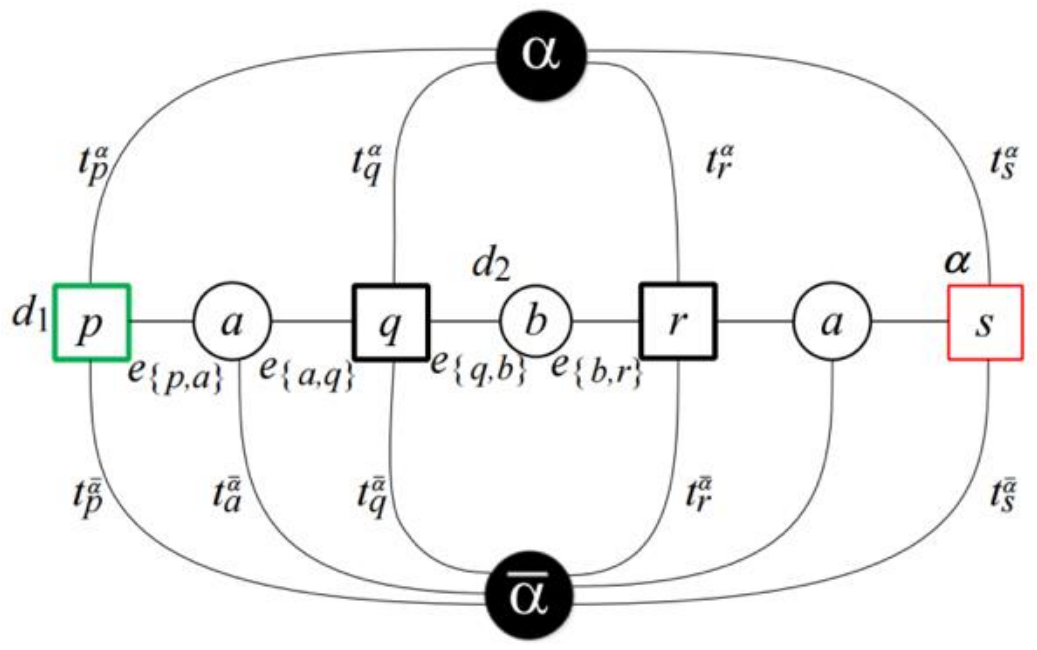

图 26　加权图的例子(以一维为例)：代表三角形的面片节点集合为$\{p, q, r, s\}$；$p$ 的标签为 $d_1$；$q$ 和 $r$ 的标签为 $d_2$；$s$ 的标签为 $\alpha$

遵循这一思想的第一个方法是由 ALEMANNO 等[90]提出的曲面分片方法。然而该方法未考虑曲面的装卡环节，需要用户人为指定各个曲面片的装卡加工方向。在此基础上，文献[88]总结相关加工特性约束(约束 1，2，4)，给出可用于减材制造加工的曲面分片方法，并通过最小化曲面片边界长度来间接满足约束 4。之后，该方法通过局部满足误差的微小变形，来保证各曲面片满足约束 3。

FANNI 等[91]提出了一种基于多边形立方体映射(Polycube Mappings)的曲面分片方法，并分析了增材制造和减材制造分片的可加工性，同时考虑了 3 轴和 4 轴铣削。由于该方法不考虑与减材制造相关的特定约束，所以最终分片结果在实践中证明仅对增材制造是可行的。MUNTONI 等[92]提出的方法通过将曲面分片这个离散的约束优化问题重新定义为一个连续函数的无约束优化，从而计算得到一组可能重叠的高度域曲面片，并且其共同覆盖了模型表面。之后，将无重叠的分片计算问题转换为连通图上的排序问题，并通过循环消除和拓扑排序的组合来求解。该组合算法计算得到了一组紧凑的高度域，在用户给定的误差范围内共同描述输入模型。YANG 等[93]为了提高加工效率，采用了一种更通用的双高度域(double height field, DHF)切片方法来取代单一高度域分片。各高度域由双曲面片组成：首先铣削一侧的曲面，然后使用适当的夹具翻转工件，铣削剩余的第二个曲面片。该方法使用了一种高效的由粗到细的分片过程：首先将输入曲面划分为最优数量的满足特定切削方向下双高度域约束的区块，然后将每个区块切割成满足体积约束的 DHF 切片。通过最小化 DHF 切片数量和最大化切片高度，该方法可以有效提升加工效率。然而，上述 3 种方法将待加工的模型进行了体素分解，并逐次使用 3 轴机床加工各个曲面块。之后，用户必须在加工完成后组装部件，以获得最终的 3D 模型，导致工件上存在可见的接缝，降低了加工质量，如图 27 所示。

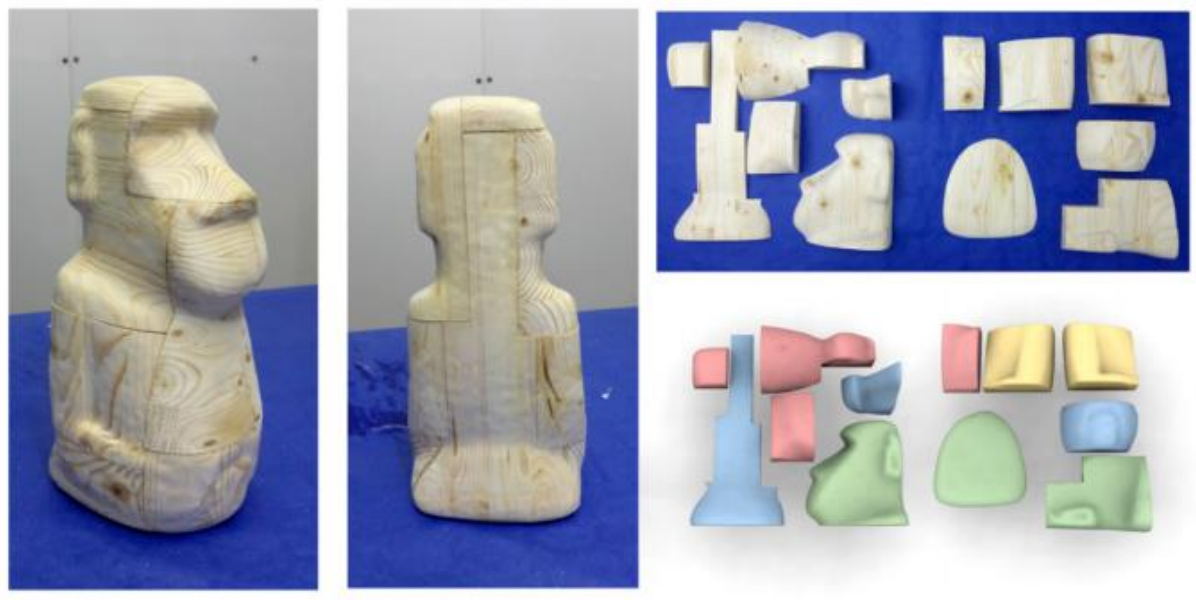

图 27　文献[92]方法的加工结果：可以看出，拼接后的模型存在明显缝隙(左图：加工结果的侧视图；中图：加工结果的背视图；右上图：加工结果的各个曲面片；右下图：曲面分片的模拟结果)

为了避免拼接过程中产生的不必要的缝隙，文献[89]提出了基于 3+2 轴机床的曲面分片算法，并提出了可达性约束(约束 3)，从而有效避免曲面片与刀具发成碰撞。该方法使用了 3+2 轴机床进行分片加工，即在加工时，2 个旋转轴先将切削刀具定位到编程位置，再由进给轴 $X$，$Y$，$Z$ 轴进行加工。在加工期间 2 个旋转轴保持不变，之后重复上述过程，直到工件被加工完成，从而有效避免 3 轴机床加工后的拼接缝隙问题。该方法首先基于约束 1，2，4 生成高度域，其次基于约束 3 生成可达性区域。之后将高度域与可达性区域利用最小覆盖集合问题(set-cover problem)相结合得到满足加工约束的曲面片。ZHAO 等[89]的方法首次保证了在多轴减材制造领域中分片算法的可行性。类似地，NUVOLI 等[94]提出了针对无边界模型的 4 轴机床曲面分片算法。此类型的机床设置扩展了经典的三轴数控加工，增加了一个围绕固定轴旋转物体的额外自由度。该方法首先根据整体加工制造精度约束确定旋转轴；其次，在旋转轴的 2 个端点处确定相应的高度域，用于装卡；最后，将网格的剩余部分分割成一组高度域，其加工方向与旋转轴垂直，并且分割时考虑了约束 4、约束 5 以及拟合逼近质量。此外，分片过程充分考虑了对象的几何特征，以及显著性信息。最终输出一组网格曲面片，可由现成的软件处理，用于生成 3 轴刀具路径。基于高度域的曲面分片方法在算法效率以及加工结果方面都有很好的表现。现有该类型算法

的基础框架均基于文献[88]的思路进行开发。然而，该算法理论是不完备的，如能量函数中的双变量项不是度量的。而根据 BOYKOV 等[95]的理论，这可能会导致最小割的权值(代价)不等于能量函数的最小值，从而导致算法的分片结果与能量函数的不对应。

## 5.2　CNC：针对加工路径的一体化插补算法

目前针对加工路径(离散 G 代码或者参数曲线)的传统进给速度规划及插补算法主要分为 2 个流程：后端：对加工路径进行特定指标下的最优速度规划(当加工路径为参数曲线时，多为离线步骤)；前端：使用后端得到的速度曲线进行实时插补，如图 28 所示。传统流程能够很好地胜任加工任务，但是存在步骤冗余繁琐、算法编程复杂、插补环节容易出现数值误差(这是由于在插补环节，一般使用高阶多项式展开，一般为二阶泰勒展开，去逼近原速度函数，而在此过程中，如果阶数太低，逼近函数与原函数之间容易产生误差)等诸多问题。因此，相关学者分别从理论驱动和数据驱动 2 个方面提出了针对加工路径的一体化插补算法，从而有效避免上述问题地产生。下面，本文将逐次介绍两类方法。

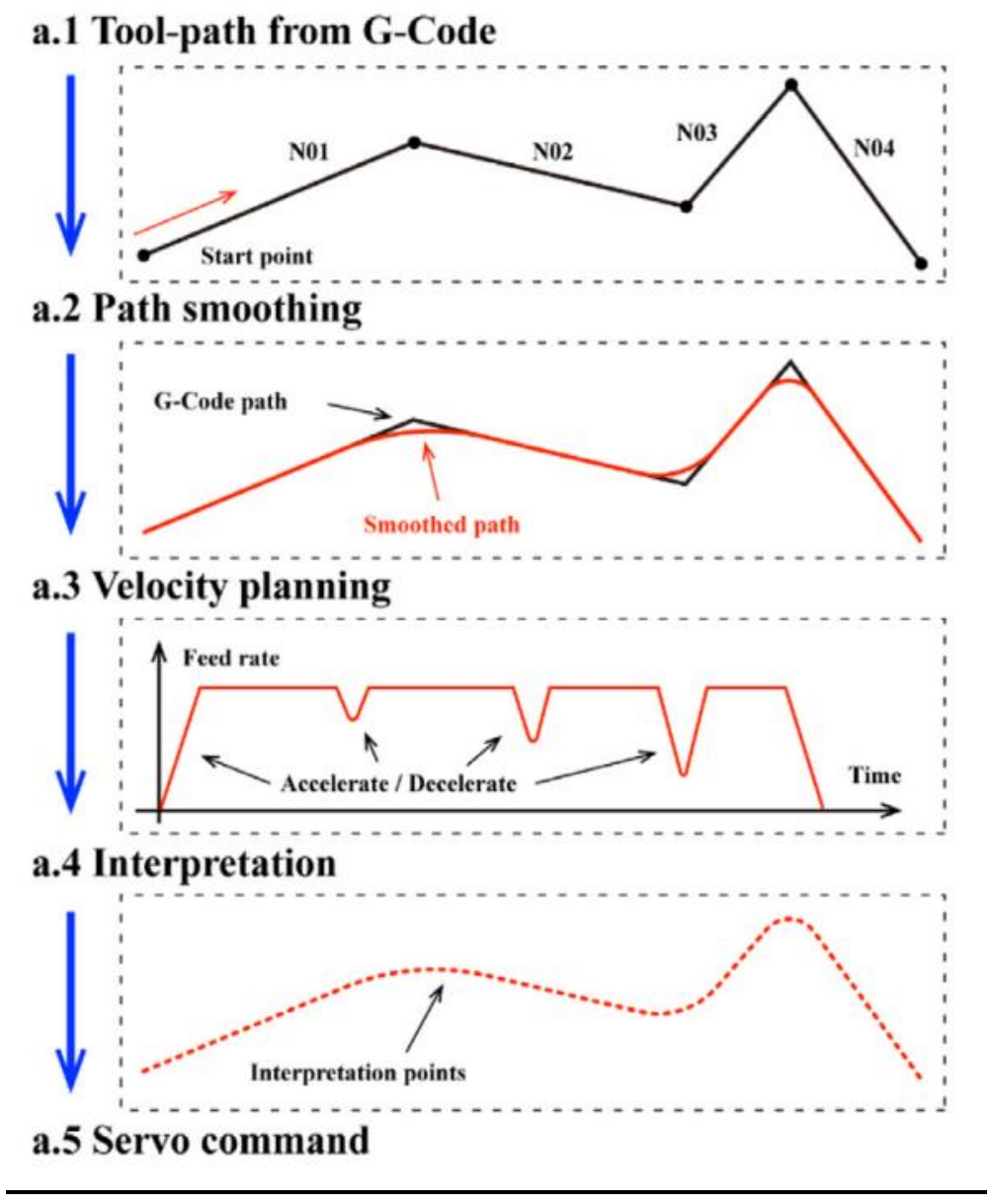

**图 28**　传统进给速度规划及插补算法流程[96]((a.1)加工路径的 G 代码；(a.2)路径光滑化；(a.3)进给速度规划；(a.4)实时插补)

### 5.2.1　理论驱动的一体化插补算法

在数控加工中，加工路径常用 G01 代码(连续的折线段)表示，直接插补会导致加工效率低和机器振荡。为避免这一问题，常采用以下 2 种方法：一是局部光滑方法，对折线路径的每个拐角做光滑过渡处理；二是整体光滑方法，将整条路径拟合为参数曲线，然后在参数曲线上规划速度。上方法均需要路径光滑和速度规划 2 个步骤，计算复杂。

YANG 等[97]提出一种新的理论驱动的一体化插补算法，即时间样条拟合方法。该方法将路径光滑化和进给速度规划合为一步完成，得到的拟合曲线以时间为参数，既包含位置信息，满足高精度的误差约束，又包含速度信息，满足运动约束，且是“bang-bang”控制的，充分利用了机床各轴的运动性能，并且可以直接用来插补。

LIN 等[98]将时间样条推广到三维情况，并引入了机床的速度、加速度和加加速度约束。为了实现高精度的轮廓误差控制，给出了折线与曲线间精确的 Hausdorff 距离计算公式。而申立勇等[99]则将时间样条曲线方法扩展到五轴数控加工中，使用五维三次 B 样条曲线去拟合机床坐标系的 G 代码。曲线由控制点和节点向量确定，曲线的 1~3 阶导数即各轴的速度、加速度和加加速度。由机床坐标系和工件坐标系的变换关系，可得到拟合曲线在工件坐标系的参数表示。因此可通过控制工件坐标系下拟合曲线与数据点的误差和机床坐标系下拟合曲线的速度、加速度以及加加速度，以加工时间最优为目标，优化样条曲线的控制点和节点向量，得到一条充分利用机床运动性能和高精度误差控制的路径。综合上述方法，袁春明等[100]提出了时间样条曲面的概念，即利用一个参数方向来描述加工路径(对原有曲面进行重新参数化，使得等参数对应于等残高)，而每条加工路径则用时间样条曲线来描述，从而将 CAD，CAM，CNC 统一起来，实现对自由曲面的一体化加工。

综上，时间样条方法可用于三轴和五轴联动加工，充分考虑加工过程中的相关约束，采用整体优化方法，求解最优的拟合路径和速度曲线，最终可由拟合曲线和插补周期直接得到插补点，是实现一体化的高精、高效加工新模式。

### 5.2.2　数据驱动的一体化插补算法

随着近几年人工智能的快速发展以及在工业领域中的广泛应用，越来越多的学者投入到其中进行研究，从而对多个行业产生了深远影响。因此，基于数据驱动的一系列研究快速崛起。基于此，文

献[96]将强化学习引入到速度规划中，提出一种基于神经网络的插补光滑算法(neural-network numerical control, NNC)来处理G代码加工路径，以实现CNC的高效一体化加工。图29展示了NNC的算法框架。可以看出，NNC直接利用训练好的神经网络来输出每个周期的插补点，从而省去中间很多的计算步骤，实时性更强。所以该方法可以分为训练过程和运行过程2个部分。而对于训练部分，算法利用强化学习来训练神经网络，因为输入的原始数据只有加工路径，属于无标签数据，所以需要使用无监督学习方式。

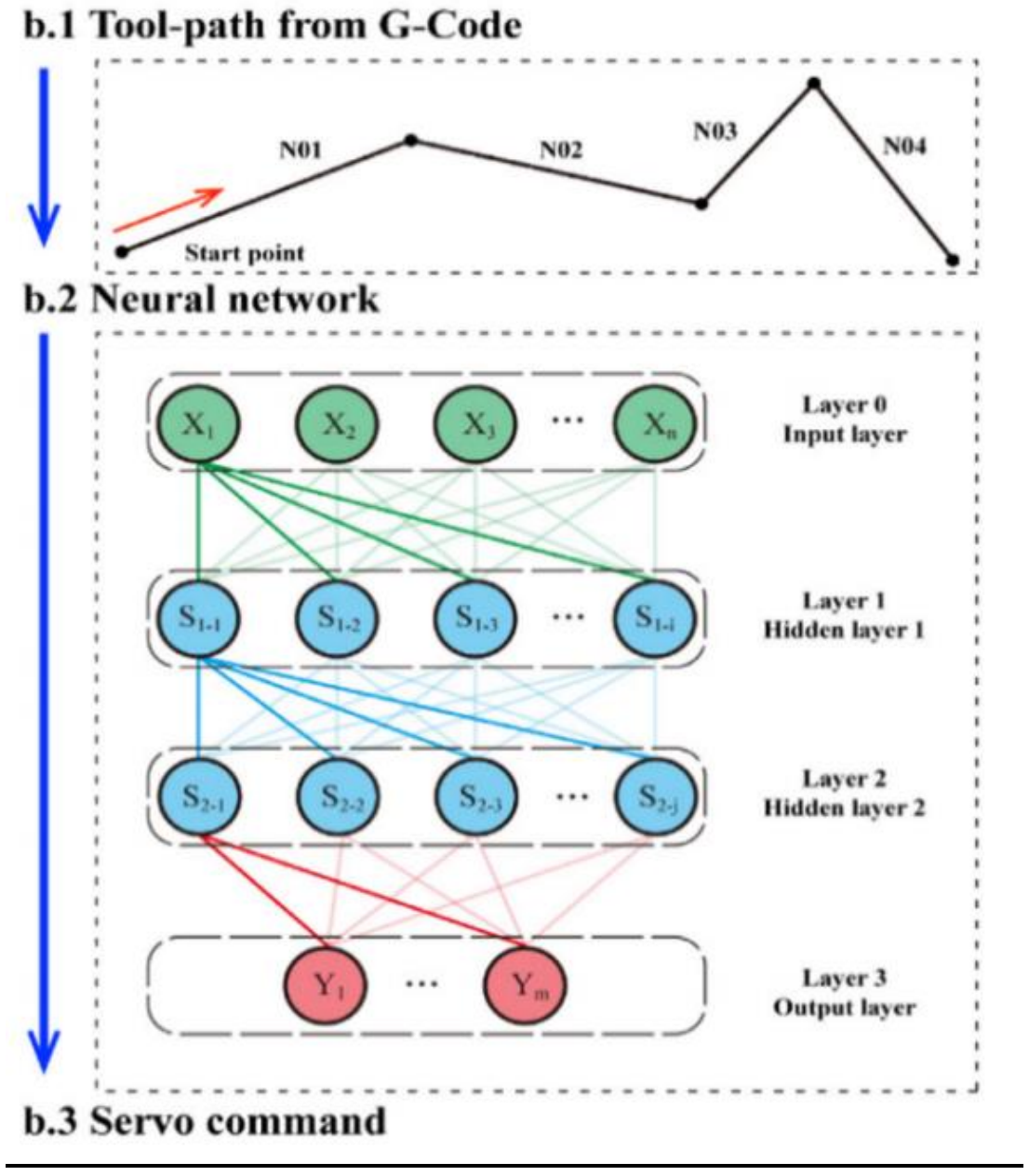

图29　NNC算法框架[96]

强化学习是一种试错方法，其目标是让智能体(Agent)在特定环境中能够采取回报最大化的行动，做出当前最优决策。其组成部分主要包括：智能体、行动(Action)、环境(Environment)、状态(State)和奖励(Reward)。其中，环境是能够将当前状态下采取的行动转换成下一个状态和奖励的函数；而智能体(一般为神经网络)则可以将新的状态和奖励转换成下一步行动。可设定智能体的表示形式，但是无法知悉环境的显式表达。所以，环境是一个只能看到输入输出的黑盒子，而强化学习则相当于尝试用智能体逼近环境函数，这样就能根据其做出最优决策。

此类型算法近几年展示出了巨大的潜力，其优势在于可以高效实时地进行一体化加工插补，但是其未充分考虑运动学约束，导致路径的加速度曲线发生超界，引起机床震颤，从而降低加工质量与效率，且该方法只能处理简单的三轴加工路径。所以此类方法未来的一个发展趋势是，构建具有微分结构的神经网络，使得生成的速度曲线满足各项运动学约束，并且考虑五轴G代码的路径加工。

## 6　未来展望

目前我国处于制造业突破瓶颈的关键时期，针对自由曲面的高效高精数控加工理论与核心技术的突破能有效提升我国工业核心竞争力。随着数控加工技术的快速发展，人们对其加工效率以及加工质量的要求也越来越严格。而当前的数控加工流程是自由曲面(CAD)—加工路径(CAM)—插补点(CNC)，即面-线-点。随着流程的逐步递进，工件在CAD阶段的一些重要几何特征(如高阶曲率信息)会逐渐丢失。目前对于比较复杂的自由曲面，人们往往以三角网格表示形式作为CAD的输出。三角网格结构简单，应用成熟，但是使用其表示的自由曲面会丢失高阶几何连续性，而这些信息正是CAM阶段所需要的。然而随着几何特征的丢失，只能依据现有的三角网格拓扑性质，通过各种拟合手段，来近似逼近原曲面的几何信息。因此会造成较大的工件精度误差，并且带来大量不必要的计算负担，进而严重降低加工效率。同时，在CAD，CAM和CNC阶段，均需要通过求解一系列优化问题，令求解对象满足各项加工特性约束，从而使得整个流程结构冗余，操作繁琐。

所以数控加工行业未来的一个重要发展趋势是将CAD，CAM和CNC深度融合，建立几何设计[101]、性能仿真[102]、设计与工艺优化[103,104]、逆向工程[105]、等流程一体化的系统，从而打破传统生产模式中“设计-规划-加工”中存在的壁垒。也就是说，在设计与优化曲面表示时，可通过加入加工特性约束以及机床性能约束，求解多目标优化问题，并得到面向CAM的可加工样条曲面。该曲面可以生成一簇等参线，两两之间满足CAM中的一系列加工特性约束，并且各等参线以加工时间为参数，在满足运动学约束和几何误差约束的前提下，达到时间最优，进而可以直接输入到机床进行插补加工。同时，相关研究可以从理论驱动和数据驱动2个方面展开，实现“数据驱动理论，理论指导数据”的完整研究链条。

参考文献 (References)

[1] CHIOU C J, LEE Y S. A machining potential field approach

to tool path generation for multi-axis sculptured surface machining[J]. Computer-Aided Design, 2002, 34(5): 357-371.
[2] LASEMI A, XUE D Y, GU P H. Recent development in CNC machining of freeform surfaces: a state-of-the-art review[J]. Computer-Aided Design, 2010, 42(7): 641-654.
[3] LONEY G C, OZSOY T M. NC machining of free form surfaces[J]. Computer-Aided Design, 1987, 19(2): 85-90.
[4] ELBER G, COHEN E. Toolpath generation for freeform surface models[J]. Computer-Aided Design, 1994, 26(6): 490-496.
[5] HE W, LEI M, BIN H Z. Iso-parametric CNC tool path optimization based on adaptive grid generation[J]. The International Journal of Advanced Manufacturing Technology, 2009, 41(5): 538-548.
[6] SUN Y W, GUO D M, JIA Z Y, et al. Iso-parametric tool path generation from triangular meshes for free-form surface machining[J]. The International Journal of Advanced Manufacturing Technology, 2006, 28(7): 721-726.
[7] ZOU Q, ZHAO J B. Iso-parametric tool-path planning for point clouds[J]. Computer-Aided Design, 2013, 45(11): 1459-1468.
[8] DING S, MANNAN M A, POO A N, et al. Adaptive iso-planar tool path generation for machining of free-form surfaces[J]. Computer-Aided Design, 2003, 35(2): 141-153.
[9] KIM B H, CHOI B K. Guide surface based tool path generation in 3-axis milling: an extension of the guide plane method[J]. Computer-Aided Design, 2000, 32(3): 191-199.
[10] YANG D O, FENG H Y. Machining triangular mesh surfaces via mesh offset based tool paths[J]. Computer-Aided Design and Applications, 2008, 5(1-4): 254-265.
[11] FENG H Y, TENG Z J. Iso-planar piecewise linear NC tool path generation from discrete measured data points[J]. Computer-Aided Design, 2005, 37(1): 55-64.
[12] HU P C, CHEN L F, TANG K. Efficiency-optimal iso-planar tool path generation for five-axis finishing machining of freeform surfaces[J]. Computer-Aided Design, 2017, 83: 33-50.
[13] KUMAZAWA G H, FENG H Y, BARAKCHI FARD M J. Preferred feed direction field: a new tool path generation method for efficient sculptured surface machining[J]. Computer-Aided Design, 2015, 67-68: 1-12.
[14] FENG H Y, LI H W. Constant scallop-height tool path generation for three-axis sculptured surface machining[J]. Computer-Aided Design, 2002, 34(9): 647-654.
[15] TOURNIER C, DUC E. A surface based approach for constant scallop HeightTool-path generation[J]. The International Journal of Advanced Manufacturing Technology, 2002, 19(5): 318-324.
[16] WEN H, GAO J, XIANG K, et al. Cutter location path generation through an improved algorithm for machining triangular mesh[J]. Computer-Aided Design, 2017, 87: 29-40.
[17] LEE E. Contour offset approach to spiral toolpath generation with constant scallop height[J]. Computer-Aided Design, 2003, 35(6): 511-518.
[18] XU J T, JI Y K, SUN Y W, et al. Spiral tool path generation method on mesh surfaces guided by radial curves[J]. Journal of Manufacturing Science and Engineering, 2018, 140(7): 071016.
[19] TRICOCHE X, ZHENG X Q, PANG A. Visualizing the topology of symmetric, second-order, time-varying two-dimensional tensor fields[M]//Mathematics and Visualization. Berlin: Springer Berlin, 2006: 225-240.
[20] MOODLEAH S, MAKHANOV S S. 5-axis machining using a curvilinear tool path aligned with the direction of the maximum removal rate[J]. The International Journal of Advanced Manufacturing Technology, 2015, 80(1): 65-90.
[21] ZHANG K, TANG K. An efficient greedy strategy for five-axis tool path generation on dense triangular mesh[J]. The International Journal of Advanced Manufacturing Technology, 2014, 74(9): 1539-1550.
[22] XU K, TANG K. Five-axis tool path and feed rate optimization based on the cutting force–area quotient potential field[J]. The International Journal of Advanced Manufacturing Technology, 2014, 75(9): 1661-1679.
[23] SU C, JIANG X, HUO G Y, et al. Initial tool path selection of the iso-scallop method based on offset similarity analysis for global preferred feed directions matching[J]. The International Journal of Advanced Manufacturing Technology, 2020, 106(7): 2675-2687.
[24] MA J W, LU X, LI G L, et al. Toolpath topology design based on vector field of tool feeding direction in sub-regional processing for complex curved surface[J]. Journal of Manufacturing Processes, 2020, 52: 44-57.
[25] ZOU Q. Length-optimal tool path planning for freeform surfaces with preferred feed directions based on Poisson formulation[J]. Computer-Aided Design, 2021, 139: 103072.
[26] PATELOUP V, DUC E, RAY P. Bspline approximation of circle arc and straight line for pocket machining[J]. Computer-Aided Design, 2010, 42(9): 817-827.
[27] ZOU Q, ZHANG J, DENG B, ZHAO, J. Iso-level tool path

planning for free-form surfaces [J]. International Journal of Production Research, Computer-Aided Design, 2014, 53: 117-125.
[28] ZOU Q. Robust and efficient tool path generation for machining low-quality triangular mesh surfaces[J]. International Journal of Production Research, 2021, 59(24): 7457-7467.
[29] AGRAWAL R K, PRATIHAR D K, ROY CHOUDHURY A. Optimization of CNC isoscallop free form surface machining using a genetic algorithm[J]. International Journal of Machine Tools and Manufacture, 2006, 46(7-8): 811-819.
[30] CRIPPS R J, CROSS B, HUNT M, et al. Singularities in five-axis machining: cause, effect and avoidance[J]. International Journal of Machine Tools and Manufacture, 2017, 116: 40-51.
[31] PLAKHOTNIK D, LAUWERS B. Graph-based optimization of five-axis machine tool movements by varying tool orientation[J]. The International Journal of Advanced Manufacturing Technology, 2014, 74(1): 307-318.
[32] BI Q Z, SHI J, WANG Y H, et al. Analytical curvature-continuous dual-Bézier corner transition for five-axis linear tool path[J]. International Journal of Machine Tools and Manufacture, 2015, 91: 96-108.
[33] TULSYAN S, ALTINTAS Y. Local toolpath smoothing for five-axis machine tools[J]. International Journal of Machine Tools and Manufacture, 2015, 96: 15-26.
[34] YANG J X, YUEN A. An analytical local corner smoothing algorithm for five-axis CNC machining[J]. International Journal of Machine Tools and Manufacture, 2017, 123: 22-35.
[35] HUANG X Y, ZHAO F, TAO T, et al. A novel local smoothing method for five-axis machining with time-synchronization feedrate scheduling[J]. IEEE Access, 2020, 8: 89185-89204.
[36] BEUDAERT X, LAVERNHE S, TOURNIER C. 5-axis local corner rounding of linear tool path discontinuities[J]. International Journal of Machine Tools and Manufacture, 2013, 73: 9-16.
[37] SHI J, BI Q Z, ZHU L M, et al. Corner rounding of linear five-axis tool path by dual PH curves blending[J]. International Journal of Machine Tools and Manufacture, 2015, 88: 223-236.
[38] XU F Y, SUN Y W. A circumscribed corner rounding method based on double cubic B-splines for a five-axis linear tool path[J]. The International Journal of Advanced Manufacturing Technology, 2018, 94(1): 451-462.
[39] MORISHIGE K, KASE K, TAKEUCHI Y. Collision-free tool path generation using 2-dimensional C-space for 5-axis control machining[J]. The International Journal of Advanced Manufacturing Technology, 1997, 13(6): 393-400.
[40] JUN C S, CHA K, LEE Y S. Optimizing tool orientations for 5-axis machining by configuration-space search method[J]. Computer-Aided Design, 2003, 35(6): 549-566.
[41] WANG N, TANG K. Automatic generation of gouge-free and angular-velocity-compliant five-axis toolpath[J]. Computer-Aided Design, 2007, 39(10): 841-852.
[42] LU J, CHEATHAM R, JENSEN C G, et al. A three-dimensional configuration-space method for 5-axis tessellated surface machining[J]. International Journal of Computer Integrated Manufacturing, 2008, 21(5): 550-568.
[43] LIN Z W, FU J Z, SHEN H Y, et al. Non-singular tool path planning by translating tool orientations in C-space[J]. The International Journal of Advanced Manufacturing Technology, 2014, 71(9): 1835-1848.
[44] HO M C, HWANG Y R, HU C H. Five-axis tool orientation smoothing using quaternion interpolation algorithm[J]. International Journal of Machine Tools and Manufacture, 2003, 43(12): 1259-1267.
[45] CASTAGNETTI C, DUC E, RAY P. The Domain of Admissible Orientation concept: a new method for five-axis tool path optimisation[J]. Computer-Aided Design, 2008, 40(9): 938-950.
[46] YE T, XIONG C H, XIONG Y L, et al. Tool orientation optimization considering second order kinematical performance of the multi-axis machine[J]. Journal of Manufacturing Science and Engineering, 2010, 132(5): 1.
[47] MI Z P, YUAN C M, MA X H, et al. Tool orientation optimization for 5-axis machining with C-space method[J]. The International Journal of Advanced Manufacturing Technology, 2017, 88(5): 1243-1255.
[48] YUAN C M, MI Z P, JIA X H, et al. Tool orientation optimization and path planning for 5-axis machining[J]. Journal of Systems Science and Complexity, 2021, 34(1): 83-106.
[49] ZHANG Y, LI Y G, XU K. Reinforcement learning–based tool orientation optimization for five-axis machining[J].The International Journal of Advanced Manufacturing Technology, 2022, 119(11-12): 7311-7326.
[50] HAN J, JIANG Y, TIAN X Q, et al. A local smoothing interpolation method for short line segments to realize continuous motion of tool axis acceleration[J]. The

International Journal of Advanced Manufacturing Technology, 2018, 95(5): 1729-1742.

[51] SUN S J, ALTINTAS Y. A G3 continuous tool path smoothing method for 5-axis CNC machining[J]. CIRP Journal of Manufacturing Science and Technology, 2021, 32: 529-549.

[52] ZHANG Y, ZHAO M Y, YE P Q, et al. A G4 continuous B-spline transition algorithm for CNC machining with jerk-smooth feedrate scheduling along linear segments[J]. Computer-Aided Design, 2019, 115: 231-243.

[53] DUAN M L, OKWUDIRE C. Minimum-time cornering for CNC machines using an optimal control method with NURBS parameterization[J]. The International Journal of Advanced Manufacturing Technology, 2016, 85(5): 1405-1418.

[54] LIN M T, LEE J C, SHEN C C, et al. Local corner smoothing with kinematic and real-time constraints for five-axis linear tool path[C]//2018 IEEE/ASME International Conference on Advanced Intelligent Mechatronics. New York: IEEE Press, 2018: 816-821.

[55] LI H X, JIANG X, HUO G Y, et al. A novel feedrate scheduling method based on Sigmoid function with chord error and kinematic constraints[J]. The International Journal of Advanced Manufacturing Technology, 2022, 119(3): 1531-1552.

[56] ZHANG L X, SUN R Y, GAO X S, et al. High speed interpolation for micro-line trajectory and adaptive real-time look-ahead scheme in CNC machining[J]. Science China Technological Sciences, 2011, 54(6): 1481-1495.

[57] CHEN Z C, KHAN M A. A new approach to generating arc length parameterized NURBS tool paths for efficient three-axis machining of smooth, accurate sculptured surfaces[J]. The International Journal of Advanced Manufacturing Technology, 2014, 70(5): 1355-1368.

[58] TULSYAN S, ALTINTAS Y. Local toolpath smoothing for five-axis machine tools[J]. International Journal of Machine Tools and Manufacture, 2015, 96: 15-26.

[59] MIN K, LEE C H, YAN C Y, et al. Six-dimensional B-spline fitting method for five-axis tool paths[J]. The International Journal of Advanced Manufacturing Technology, 2020, 107(5): 2041-2054.

[60] ZHANG M, YAN W, YUAN C M, et al. Curve fitting and optimal interpolation on CNC machines based on quadratic B-splines[J]. Science China Information Sciences, 2011, 54(7): 1407-1418.

[61] TIMAR S D, FAROUKI R T, SMITH T S, et al. Algorithms for time-optimal control of CNC machines along curved tool paths[J]. Robotics and Computer-Integrated Manufacturing, 2005, 21(1): 37-53.

[62] BIAGIOTTI L, MELCHIORRI C. FIR filters for online trajectory planning with time- and frequency-domain specifications[J]. Control Engineering Practice, 2012, 20(12): 1385-1399.

[63] TAJIMA S, SENCER B, SHAMOTO E. Accurate interpolation of machining tool-paths based on FIR filtering[J]. Precision Engineering, 2018, 52: 332-344.

[64] SENCER B, KAKINUMA Y, YAMADA Y. Linear Interpolation of machining tool-paths with robust vibration avoidance and contouring error control[J]. Precision Engineering, 2020, 66: 269-281.

[65] TAJIMA S, SENCER B. Accurate real-time interpolation of 5-axis tool-paths with local corner smoothing[J]. International Journal of Machine Tools and Manufacture, 2019, 142: 1-15.

[66] LIU Y, WAN M, QIN X B, et al. FIR filter-based continuous interpolation of G01 commands with bounded axial and tangential kinematics in industrial five-axis machine tools[J]. International Journal of Mechanical Sciences, 2020, 169: 105325.

[67] YE P Q, SHI C, YANG K M, et al. Interpolation of continuous micro line segment trajectories based on look-ahead algorithm in high-speed machining[J]. The International Journal of Advanced Manufacturing Technology, 2008, 37(9): 881-897.

[68] ERKORKMAZ K, ALTINTAS Y. High speed CNC system design. Part I: jerk limited trajectory generation and quintic spline interpolation[J]. International Journal of Machine Tools and Manufacture, 2001, 41(9): 1323-1345.

[69] JAHANPOUR J, ALIZADEH M R. A novel acc-jerk-limited NURBS interpolation enhanced with an optimized S-shaped quintic feedrate scheduling scheme[J]. The International Journal of Advanced Manufacturing Technology, 2015, 77(9): 1889-1905.

[70] JIA Z Y, SONG D N, MA J W, et al. A NURBS interpolator with constant speed at feedrate-sensitive regions under drive and contour-error constraints[J]. International Journal of Machine Tools and Manufacture, 2017, 116: 1-17.

[71] FAN W, GAO X S, YAN W, et al. Interpolation of parametric CNC machining path under confined jounce[J]. The International Journal of Advanced Manufacturing Technology, 2012, 62(5): 719-739.

[72] ZHANG L Q, DU J F. Acceleration smoothing algorithm based on jounce limited for corner motion in high-speed machining[J]. The International Journal of Advanced Manufacturing Technology, 2018, 95(1): 1487-1504.

[73] WANG Y S, YANG D S, GAI R L, et al. Design of trigonometric velocity scheduling algorithm based on pre-interpolation and look-ahead interpolation[J]. International Journal of Machine Tools and Manufacture, 2015, 96: 94-105.

[74] SUN Y W, ZHAO Y, BAO Y R, et al. A novel adaptive-feedrate interpolation method for NURBS tool path with drive constraints[J]. International Journal of Machine Tools and Manufacture, 2014, 77: 74-81.

[75] SUN Y W, CHEN M S, JIA J J, et al. Jerk-limited feedrate scheduling and optimization for five-axis machining using new piecewise linear programming approach[J]. Science China: Technological Sciences, 2019, 62(7): 1067-1081.

[76] SUN Y W, SHI Z F, XU J T. Synchronous feedrate scheduling for the dual-robot machining of complex surface parts with varying wall thickness[J]. The International Journal of Advanced Manufacturing Technology, 2022, 119(3): 2653-2667.

[77] LIU H, LIU Q, SUN P P, et al. The optimal feedrate planning on five-axis parametric tool path with geometric and kinematic constraints for CNC machine tools[J]. International Journal of Production Research, 2017, 55(13): 3715-3731.

[78] YE P Q, ZHANG Y, XIAO J X, et al. A novel feedrate planning and interpolating method for parametric toolpath in Frenet-Serret frame[J]. The International Journal of Advanced Manufacturing Technology, 2019, 101(5): 1915-1925.

[79] YUAN C M, ZHANG K, FAN W. Time-optimal interpolation for CNC machining along curved tool pathes with confined chord error[J]. Journal of Systems Science and Complexity, 2013, 26(5): 836-870.

[80] LIU Y, WAN M, QIN X B, et al. FIR filter-based continuous interpolation of G01 commands with bounded axial and tangential kinematics in industrial five-axis machine tools[J]. International Journal of Mechanical Sciences, 2020, 169: 105325.

[81] TAJIMA S, SENCER B. Real-time trajectory generation for 5-axis machine tools with singularity avoidance[J]. CIRP Annals, 2020, 69(1): 349-352.

[82] ZHANG Q, GAO X S, LI H B, et al. Minimum time corner transition algorithm with confined feedrate and axial acceleration for nc machining along linear tool path[J]. The International Journal of Advanced Manufacturing Technology, 2017, 89(1): 941-956.

[83] TAJIMA S, SENCER B. Kinematic corner smoothing for high speed machine tools[J]. International Journal of Machine Tools and Manufacture, 2016, 108: 27-43.

[84] ZHANG Y B, WANG T Y, PENG P, et al. Feedrate blending method for five-axis linear tool path under geometric and kinematic constraints[J]. International Journal of Mechanical Sciences, 2021, 195: 106262.

[85] WANG W X, HU C X, ZHOU K, et al. Local asymmetrical corner trajectory smoothing with bidirectional planning and adjusting algorithm for CNC machining[J]. Robotics and Computer-Integrated Manufacturing, 2021, 68: 102058.

[86] XIAO J L, LIU S J, LIU H T, et al. A jerk-limited heuristic feedrate scheduling method based on particle swarm optimization for a 5-DOF hybrid robot[J]. Robotics and Computer-Integrated Manufacturing, 2022, 78: 102396.

[87] LIANG F S, YAN G P, FANG F Z. Global time-optimal B-spline feedrate scheduling for a two-turret multi-axis NC machine tool based on optimization with genetic algorithm[J]. Robotics and Computer-Integrated Manufacturing, 2022, 75: 102308.

[88] HERHOLZ P, MATUSIK W, ALEXA M. Approximating free-form geometry with height fields for manufacturing[J]. Computer Graphics Forum, 2015, 34(2): 239-251.

[89] ZHAO H S, ZHANG H, XIN S Q, et al. DSCarver: decompose-and-spiral-carve for subtractive manufacturing[J]. ACM Transactions on Graphics, 2018, 37(4): 137.

[90] ALEMANNO G, CIGNONI P, PIETRONI N, et al. Interlocking pieces for printing tangible cultural heritage replicas[C]//Proceedings of the Eurographics Workshop on Graphics and Cultural Heritage. New York: ACM Press, 2014: 145-154.

[91] FANNI F A, CHERCHI G, MUNTONI A, et al. Fabrication oriented shape decomposition using polycube mapping[J]. Computers & Graphics, 2018, 77: 183-193.

[92] MUNTONI A, LIVESU M, SCATENI R, et al. Axis-aligned height-field block decomposition of 3D shapes[J]. ACM Transactions on Graphics, 2018, 37(5): 169.

[93] YANG J F, ARAUJO C, VINING N, et al. DHFSlicer: double height-field slicing for milling fixed-height materials[J]. ACM Transactions on Graphics, 2020, 39(6): 205.

[94] NUVOLI S, TOLA A, MUNTONI A, et al. Automatic surface

segmentation for seamless fabrication using 4 - axis milling machines[J]. Computer Graphics Forum, 2021, 40(2): 191-203.

[95] BOYKOV Y, VEKSLER O, ZABIH R. Fast approximate energy minimization via graph cuts[J]. IEEE Transactions on Pattern Analysis and Machine Intelligence, 2001, 23(11): 1222-1239.

[96] LI B R, ZHANG H, YE P Q, et al. Trajectory smoothing method using reinforcement learning for computer numerical control machine tools[J]. Robotics and Computer-Integrated Manufacturing, 2020, 61: 101847.

[97] YANG Z Y, SHEN L Y, YUAN C M, et al. Curve fitting and optimal interpolation for CNC machining under confined error using quadratic B-splines[J]. Computer-Aided Design, 2015, 66: 62-72.

[98] LIN F M, SHEN L Y, YUAN C M, et al. Certified space curve fitting and trajectory planning for CNC machining with cubic B-splines[J]. Computer-Aided Design, 2019, 106: 13-29.

[99] 申立勇, 袁春明, 高小山, 等. 一种应用于五轴数控机床的时间样条曲线拟合与插补方法: CN114115131A[P]. 2022-03-01.

SHEN L Y, YUAN C M, GAO X S, et al. A time-spline curve fitting and interpolation method applied to five-axis CNC machine tools: CN114115131A[P]. 2022-03-01.

[100] 袁春明, 申立勇, 高小山, 等. 一种基于 CAM 的时间样条曲面生成方法: CN114217572A[P]. 2022-03-22.

YUAN C M, SHEN L Y, GAO X S, et al. A time-spline surface generation method based on CAM: CN114217572A[P]. 2022-03-22.

[101] ZOU Q, FENG HY. Push-pull direct modeling of solid CAD models[J]. Advances in Engineering Software. 2019, 127:59-69.

[102] HUGHS TJ, COTTRELL JA, BAZILEVIS Y. Isogeometric analysis: CAD, finite elements, NURBS, exact geometry and mesh refinement[J]. Computer methods in applied mechanics and engineering. 2005, 194(39-41):4135-95.

[103] WANG G, ZOU Q, ZHAO C, LIU Y, YE X. A highly efficient approach for bi-level programming problems based on dominance determination[J]. Journal of Computing and Information Science in Engineering. 2022, 22(4):041006.

[104] ZUBAIR AF, MANSOR MS. Embedding firefly algorithm in optimization of CAPP turning machining parameters for cutting tool selections[J]. Computers & Industrial Engineering. 2019, 135:317-25.

[105] SU C, JIANG X, HUO G, ZOU Q, ZHENG Z, FENG HY. Accurate model construction of deformed aero-engine blades for remanufacturing[J]. The International Journal of Advanced Manufacturing Technology. 2020, 106:3239-51.